\documentclass[10pt,notitlepage,twoside,a4paper]{amsart}

\usepackage{amsmath,amssymb,enumerate}

\usepackage{epsfig,fancyhdr,color}

\usepackage{amssymb}
\usepackage{amsmath,amsthm}
\usepackage{latexsym}
\usepackage{amscd}
\usepackage{psfrag}
\usepackage{graphicx}
\usepackage[latin1]{inputenc}
\usepackage[all]{xy}
\usepackage[mathcal]{eucal}

\newtheorem{theorem}{\rm\bf Theorem}[section]
\newtheorem{proposition}[theorem]{\rm\bf Proposition}
\newtheorem{lemma}[theorem]{\rm\bf Lemma}
\newtheorem{corollary}[theorem]{\rm\bf Corollary}
\newtheorem*{theorem 1}{\rm\bf Proposition 1}
\newtheorem*{theorem 2}{\rm\bf Proposition 2}

\theoremstyle{definition}
\newtheorem{definition}[theorem]{\rm\bf Definition}

\theoremstyle{remark}
\newtheorem{remark}[theorem]{\rm\bf Remark}

\newtheorem{example}[theorem]{\rm\bf Example}

\newtheorem{question}[theorem]{\rm\bf Question}

\def\interieur#1{\mathord{\mathop{\kern 0pt #1}\limits^\circ}}

\title[ Teichm\"uller spaces of surfaces of infinite type]{Some metrics on Teichm\"uller spaces of surfaces of infinite type}

\author{Lixin Liu}
\address{L. Liu, Department of Mathematics, Zhongshan University, 510275, Guangzhou, P. R. China}
\email{mcsllx@mail.sysu.edu.cn}

\author{Athanase Papadopoulos}
\address{A. Papadopoulos, Institut de Recherche Math{\'e}matique Avanc\'ee,
Universit{\'e} de Strasbourg and CNRS,
7 rue Ren\'e Descartes,
 67084 Strasbourg Cedex - France} \email{papadopoulos@math.u-strasbg.fr}
\date{\today}

\begin{document}

\begin{abstract}   Unlike the case of surfaces of topologically finite type, there are several different Teichm\"uller spaces that are associated to a surface of topological infinite type. These Teichm\"uller spaces first depend (set-theoretically) on whether we work in the hyperbolic category or in the conformal category. They also depend, given the
  choice of a point of view (hyperbolic or conformal),  on the choice of a distance function on Teichm\"uller space. Examples of distance functions that appear naturally in the hyperbolic setting are the length spectrum distance and the bi-Lipschitz distance, and there are other useful distance functions. The Teichm\"uller spaces  also depend on the choice of a basepoint. The aim of this paper is to present some examples, results and questions on the Teichm\"uller theory of surfaces of infinite topological type that do not appear in the setting the Teichm\"uller theory of surfaces of finite type. In particular, we point out relations and differences between the various Teichm\"uller spaces associated to a given surface of topological infinite type.

\bigskip

\medskip


\medskip
\noindent Acknowledgement:   L. Liu was partially supported by the NSF of China.

\end{abstract}
\maketitle

\section{Introduction}\label{intro}

Let $X$ be a connected oriented surface of infinite topological type. Unlike the case of a surface of finite type, there are several Teichm\"uller spaces associated to $X$. Such a Teichm\"uller space can be defined as a space of equivalence classes of marked hyperbolic surfaces homeomorphic to $X$, the marking being a homeomorphism between a base hyperbolic structure on $X$ and a hyperbolic surface homeomorphic to $X$. To define a distance function (and then, a topology) on Teichm\"uller space, we ask that this marking be either length spectrum bounded or a bi-Lipschitz homeomorphism. Depending on such a choice, the Teichm\"uller spaces that we obtain are setwise different.
Taking a different point of view, a Teichm\"uller space associated to $X$ can also be defined as a space of equivalence classes of marked Riemann surfaces homeomorphic to $X$, the marking being here a quasiconformal homeomorphism between a base Riemann surface structure on $X$ and a Riemann surface homeomorphic to $X$. We shall be more precise below. Unlike the case of surfaces of finite type, the various definitions of Teichm\"uller spaces are not equivalent and the resulting Teichm\"uller space strongly depends on the point of view (hyperbolic or conformal), on the distance function we choose on the set of equivalence classes of marked surfaces, and on the choice of a base surface.
The aim of this paper is to review some of the Teichm\"uller spaces associated to $X$, precisely when $X$ has infinite type, and to point out some differences, some relations and some analogies between these various Teichm\"uller spaces.

                    The plan of this paper is as follows.

          In Section \ref{s1}, we study the length spectrum Teichm\"uller spaces associated to a surface of infinite type.  This  study involves the consideration of the length spectrum distance between hyperbolic surfaces.

 In Section \ref{s22}, we study the  bi-Lipschitz Teichm\"uller spaces associated to  a hyperbolic surface of infinite type, and this involves  the consideration of the bi-Lipschitz distance between hyperbolic surfaces.

                In Section \ref{s11}, we consider the quasiconformal point of view. We define the quasiconformal  Teichm\"uller space of a Riemann surface of infinite type and the quasiconformal metric on this space. This space is related to the bi-Lipschitz Teichm\"uller space. We give a sufficient condition under which the three Teichm\"uller spaces coincide as sets.

Each section contains basic definitions and results, and some open questions.

We now fix some notation which will be used throughout the paper.

We consider a connected orientable surface $X$ of genus $g\geq 0$, with $p\geq 0$ punctures and $n\geq 0$ boundary components, and we assume that $g+n+p$ is countably infinite.
We note right away that even in the case of surfaces with no punctures and no boundary components, two surfaces of countably infinite genus are not necessarily homeomorphic. (For instance, the surfaces represented in Figures \ref{infinite}, \ref{bi-infinite} and \ref{tri-infinite} below are pairwise non-homeomorphic.) This contrasts with the case of surfaces of finite type in which the triple (genus, number of punctures, number of boundary components) determines the surface up to homeomorphism. The classification of noncompact surfaces has a long history, and for this subject we refer the reader to the work of Ker\'ekj\'art\'o \cite{surfaces2}, to the later paper by Richards \cite{surfaces3}, and to recent update by Prishlyak and Mischenko \cite{surfaces4}. In particular, there exist invariants  for surfaces of infinite type. (We note that these invariants are, unlike those of  surfaces of finite type, not discrete.)

Using Zorn's Lemma, on any topological surface $X$ as above, we can find a countable family of disjoint simple closed curves $C_i$, $i=1,2,\ldots$ such that each connected component of $S-\cup_{i=1}^{\infty} C_i$ is homeomorphic to a pair of pants, that is, a sphere with three holes. Here, a hole is either a puncture (that is, a point removed from the surface) or an open disk removed. In the latter case the surface has a boundary component at the hole.

 All the hyperbolic structures that we consider on $X$ will be complete, such that each boundary component of this surface a closed geodesic and each puncture of this surface has a neighborhood isometric to a cusp, that is, to the quotient of a region of the form $\{(x,y)\ \vert x>a\}$ in the upper half-space model of the hyperbolic plane $\mathbb{H}^2$, with $a>0$,  by the transformation $z\mapsto z+1$. The neighborhood of a cusp is biholomorphically equivalent to a punctured disk in $\mathbb{C}$.

For future reference, we reformulate our requirements as follows:

\medskip
\noindent ($\star$)  The surface $X$  admits a pants decomposition $\mathcal{P}$ with a countable number of pair of pants. The hyperbolic structures that we consider on $X$ are all metrically complete, and they have the property that if we replace each boundary component of a pair of pants in the decomposition $\mathcal{P}$ by the closed geodesic in its homotopy class, then each pair of pants becomes a sphere with three holes, a hole being either a boundary component which is a closed geodesic (and this can hold when the hole is either a boundary component of $X$ or a closed curve in the pants decomposition $\mathcal{P}$) or a cusp (and this holds when the hole is a puncture of $X$). In the case where the surface is just a Riemann surface (i.e. not necessarily equipped with a hyperbolic metric), then, we assume that the unique complete hyperbolic metric  on $X\setminus \partial X$ that uniformizes the complex structure on that surface, after cutting the flares if such cylindrical ends exist, satisfies the above conditions.

\medskip
We point out right away that the completeness of a hyperbolic surface satisfying Property  ($\star$) is not a redundant property. This contrasts with the case of hyperbolic surfaces of finite type with geodesic boundary, where a surface is complete if and only if the neighborhood of each puncture is a cusp. Examples of non-complete infinite type hyperbolic surfaces with no punctures and no boundary components have been given by Basmajian  in \cite{Basmajian}.

It is easy (although not completely trivial) to prove that any hyperbolic surface satisfying  ($\star$) has infinite diameter.

A large part of the classical quasiconformal theory of Teichm\"uller space (as started by Teichm\"uller, and developed by Ahlfors and Bers) is valid for surfaces of infinite type. In contrast, Thurston's surface theory does not easily extend to surfaces of infinite type. For instance, on any surface of infinite type,  we can find an infinite sequence of homotopy classes of simple closed curves which are disjoint, pairwise nonhomotopic and nonhomotopic to a point, and such a sequence does not converge in any reasonable sense to a homotopy class of measured foliations.

    We shall consider ordered pairs $(f,H)$, where $H$ is a hyperbolic metric on a surface $S$ homeomorphic to $X$, and $f:X\to S$ a homeomorphism. Such a pair $(f,H)$ is called a {\it marked hyperbolic surface} (in general, with a base surface $X$ being understood).
The homeomorphism $f$ is called the {\it marking} of $H$, or of $S$. We shall also consider a hyperbolic metric on $X$ as a marked hyperbolic surface, taking as marking the identity map. Finally, we shall also talk about markings $f:H_0\to H$  between hyperbolic surfaces, when  considering the hyperbolic surface $H_0$ on $S$ as a basepoint in our space.

    A {\it simple closed curve} on $X$ is a one-dimensional submanifold of $X$ homeomorphic to a circle. We shall call a curve on $X$ {\it essential} if it is a simple closed curve which is neither homotopic to a point nor to a puncture (but it can be homotopic to a boundary component).

We let $\mathcal{S}=\mathcal{S}(X)$ be the set of homotopy classes of essential curves on $X$.

The {\it intersection number} $i(\alpha,\beta)$ of two elements $\alpha$ and $\beta$ in $\mathcal{S}$ is the minimum number of intersection points of two closed curves representing the classes $\alpha$ and $\beta$.

Given a hyperbolic structure $H$ on $X$ and given a homotopy class $\alpha$ of essential curves on $X$, we denote by $l_H(\alpha)$ the length of the closed $H$-geodesic on $X$ in the class $\alpha$.

 Finally, we note that the Teichm\"uller spaces that we consider in this paper are, in the classical terminology, ``reduced Teichm\"uller spaces". We shall recall the definitions below.

We thank Daniele Alessandrini and the referee for carefully reading the original manuscript and for their very useful remarks.

    \section{The length spectrum Teichm\"ulller spaces}\label{s1}

       Throughout this section,  $H_0$  is a fixed hyperbolic metric on the surface $X$ called the base hyperbolic structure. With this, a marking $f:X\to S$ (where $S$ is a hyperbolic surface) becomes a homeomorphism between hyperbolic surfaces. Given  a homeomorphism $f$ between two surfaces $S$  and $S'$ equipped with hyperbolic metrics $H$ and $H'$ respectively, we say that $f$ is {\it length-spectrum bounded} if the following holds:

\begin{equation}\label{eq:ls}
K(f)= \sup_{\alpha\in\mathcal{S}(H)} \left\{ \frac{l_{H'}(f(\alpha))}{l_{H}(\alpha)},\frac{l_{H}(\alpha)}{l_{H'}(f(\alpha))}\right\}<\infty.
\end{equation}

We shall call the quantity $K(f)$ that appears  in (\ref{eq:ls}) the {\it length-spectrum constant of $f$}. Note that this quantity only depends on the homotopy class of $f$.

We consider the collection of marked hyperbolic structures $(f,H)$ relative to the base structure $H_0$ and where the marking
$f:H_0\to H$ is length-spectrum bounded. Given two such marked hyperbolic structures $(f,H)$ and $(f',H')$, we write $(f,H)\sim(f',H')$ if there exists an isometry (or, equivalently, a length spectrum preserving homeomorphism, see Proposition \ref{prop:isometry-} below) $f'':H\to H'$ which is homotopic to $f'\circ  f^{-1}$. We note that all the homotopies of a surface that we consider in this paper  preserve the punctures and preserve setwise the boundary components of the surface at all times.
The relation $\sim$ is an equivalence relation on the set of length-spectrum bounded marked hyperbolic surfaces $(f,H)$, based at $(X,H_0)$.

\begin{definition} The {\it length-spectrum Teichm\"uller space} $\mathcal{T}_{ls}(H_0)$ is the space of $\sim$-equivalence classes $[f,H]$ of length-spectrum bounded marked hyperbolic surfaces $(f,H)$. The {\it basepoint} of this Teichm\"uller space is the equivalence class $[\mathrm{Id},H_0]$.
\end{definition}

 We note that the fact that we do not ask our homotopies to preserve pointwise the boundary of the surface corresponds to working with what is usually called the {\it reduced} Teichm\"uller space, instead of Teichm\"uller space. (In the latter case, the homotopies that define the equivalence relation are required to induce the identity map on each boundary component.) There is a substantial difference between the two theories. For instance, in the quasiconformal setting that we consider in Section \ref{s11} below, the Teichm\"uller space of the unit disk in $\mathbb{C}$ is,  in the non-reduced theory, infinite-dimensional (and it is called the {\it universal Teichm\"uller space}), whereas in the reduced theory, this space is reduced to a point. Of course, for surfaces that do not have boundary components, the reduced and non-reduced Teichm\"uller spaces coincide.
 We also note that in conformal geometry, dealing with boundary points of a Riemann surface generally involves working with the Poincar\'e metric and its boundary theory, that is, using hyperbolic geometry.

Since all the Teichm\"uller spaces that we use in this paper are reduced, we shall use, for simplicity, the terminology {\it Teichm\"uller space} instead of {\it reduced Teichm\"uller space}.

 The topology of $\mathcal{T}_{ls}(H_0)$ is induced by the {\it length-spectrum} metric $d_{ls}$, defined by taking the distance $d_{ls}([f,H],[f',H'])$ between two points in $\mathcal{T}_{ls}(H_0)$ represented by two marked hyperbolic surfaces $(f,H)$ and $(f',H')$ to be
        \[d_{ls}([f,H],[f',H'])=\frac{1}{2}\log K(f'\circ  f^{-1}).\]
        
        (It may be useful to recall here that the length-spectrum constants  of length-spectrum bounded homeomorphism only depends on the homotopy class of such a homeomorphism.)

The fact that the function $d_{ls}$ satisfies the properties of a metric is straightforward, except perhaps for the separation axiom, which is a consequence of the following:

            \begin{proposition} \label{prop:isometry-} Let $H$ and $H'$ be two hyperbolic metrics on a surface $X$ and suppose that $f:(X,H)\to (X,H')$ is a homeomorphism whose length-spectrum constant    $K(f)$ equals 1. Then, $f$ is isotopic to an isometry.
    \end{proposition}
            \begin{proof} From the hypothesis, it follows that for each $\alpha$ in $\mathcal{S}(X)$ we have $l_{H}(\alpha)=l_{H'}(f(\alpha))$. We show that $f$ is isotopic to an isometry.
       Consider a sequence $K_i$, $i=0,1,\ldots$ of finite type subsurfaces with boundary satisfying $K_i\subset K_{i+1}$ for all $i\geq 0$ and $\cup_{i=0}^{\infty}K_i=X$.  For each $i\geq 0$, let $L_i\subset X$ (respectively $M_i\subset X$) be the hyperbolic subsurface of $X$ isotopic to $K_i$ and having geodesic boundary, with respect to the hyperbolic metric $H$ (respectively $H'$). Then, for each $i\geq 0$, $f$ is homotopic to a map $f'_i:X\to X$ which sends homeomorphically $L_i$ to $M_i$. From a well known result for surfaces of finite type, since $f'_i$ preserves the length spectrum and since each $L_i$ is of finite type, the restriction of $f'_i$ to each $L_i$ is isotopic to an isometry $f''_i:L_i\to M_i$. Now an isometry between two surfaces of finite type is  unique in its isotopy class. Thus, for each $i\geq 1$, the restriction of $f'_i$ to $L_{i-1}$ is equal to $f''_{i-1}$. By taking the union of the maps $f''_i$, we obtain an isometry from $(X,H)$ to $(X,H')$.

        \end{proof}

        The length function $\alpha\mapsto l_H(\alpha)$ on $\mathcal{S}(X)$ associated to a hyperbolic structure $H$ on $X$ induces a length function $\alpha\mapsto l_x(\alpha)$ associated to an element $x$ of the Teichm\"uller space $\mathcal{T}_{ls}(H_0)$. The first observation is the following:
        \begin{proposition} For any base hyperbolic metric $H_0$ on $X$ and for any $\alpha$ in $\mathcal{S}$, the map $x\mapsto \log l_x(\alpha)$ defined on $\mathcal{T}_{ls}(H_0)$ is 2-Lipschitz.
        \end{proposition}
        \begin{proof} Let $x$ and $x'$ be two points in $\mathcal{T}_{ls}(H_0)$, represented by marked surfaces $(f,H)$ and $(f',H')$ respectively. Then, for any $\alpha$ in $\mathcal{S}$ and for any length-spectrum bounded homeomorphism $f'':H\to H'$ with associated length-spectrum constant $K(f'')$, we have

 \begin{equation*} 
 \begin{array}{llll}
    \vert \log l_x(\alpha)- \log l_{x'}(\alpha)\vert& =& \left\vert \log \frac{l_{H'}(\alpha)}{l_{H}(\alpha)}\right\vert\\
 & \leq & \max \left\{\log \frac{l_{H'}(\alpha)}{ l_{H}(\alpha)},\log \frac{ l_{H}(\alpha)}{l_{H'}(\alpha)}\right\}\\
 & \leq&\log K(f'').
 \end{array}
\end{equation*}

        Taking the infimum over all  length-spectrum bounded homeomorphisms $f'':H\to H'$ in the homotopy class of $f'\circ f^{-1}$, we obtain
        \[\vert \log l_x(\alpha)-\log l_{x'}(\alpha)\vert\leq 2d_{ls}(x,x').\]
        \end{proof}

        The next observation is the existence of pairs of marked hyperbolic structures on $X$ which are not related by any length-spectrum bounded homeomorphism. This will show that two Teichm\"uller spaces $\mathcal{T}_{ls}(H_0)$ and $\mathcal{T}_{ls}(H_1)$ based at two different hyperbolic structures $H_0$ and $H_1$ on $X$ are in general different as sets.

        \bigskip
  \begin{figure}[!hbp]
\centering
\psfrag{a}{\small $C_1$}
\psfrag{b}{\small $C_2$}
\psfrag{c}{\small $C_3$}
\psfrag{d}{\small $\ldots$}
\includegraphics[width=.65\linewidth]{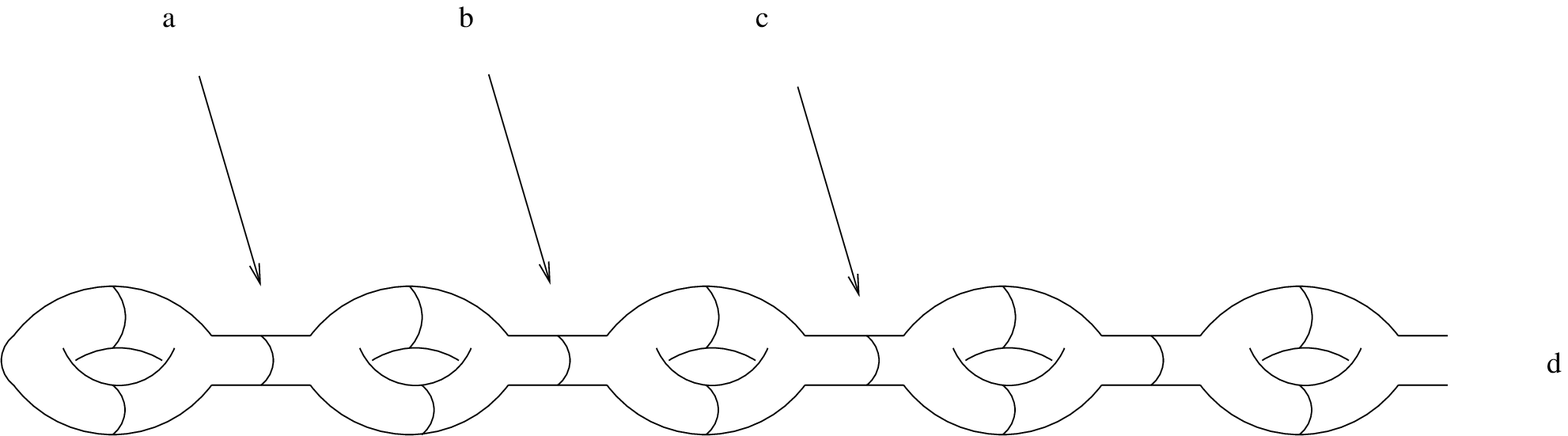}
\caption{\small {}}
\label{infinite}
\end{figure}
\bigskip

    \begin{example}\label{ex:ls}
    Consider the surface drawn in Figure \ref{infinite}, where $C_1,C_2,\ldots$ are the curves in the homotopy classes represented. We equip this surface with two
      hyperbolic metrics $H_0$ and $H_1$ such that for all $i=1,2,\ldots$,   $l_{H_{0}}(C_i)=1$ and $l_{H_{1}}(C_i)=1/i$. The metrics $H_0$ and $H_1$ have the further property that any closed ball of radius 1 on the surface is contained in a finite number of pairs of pants (of the given decomposition), and therefore it is compact. Thus, by the theorem of Hopf-Rinow, the two metrics are complete.  It is clear that the hyperbolic structures $H_0$ and $H_1$ are not length-spectrum equivalent, and therefore $H_1\not\in\mathcal{T}_{ls}(H_0)$. In particular, $\mathcal{T}_{ls}(H_0)\not=\mathcal{T}_{ls}(H_1)$.
    \end{example}

\bigskip

  \begin{figure}[!hbp]
\centering
\psfrag{a}{\small $C_{-1}$}
\psfrag{b}{\small $C_0$}
\psfrag{c}{\small $C_1$}
\psfrag{d}{\small $C_{2}$}
 \psfrag{e}{\small $\ldots$}
\psfrag{f}{\small $\ldots$}
\includegraphics[width=.65\linewidth]{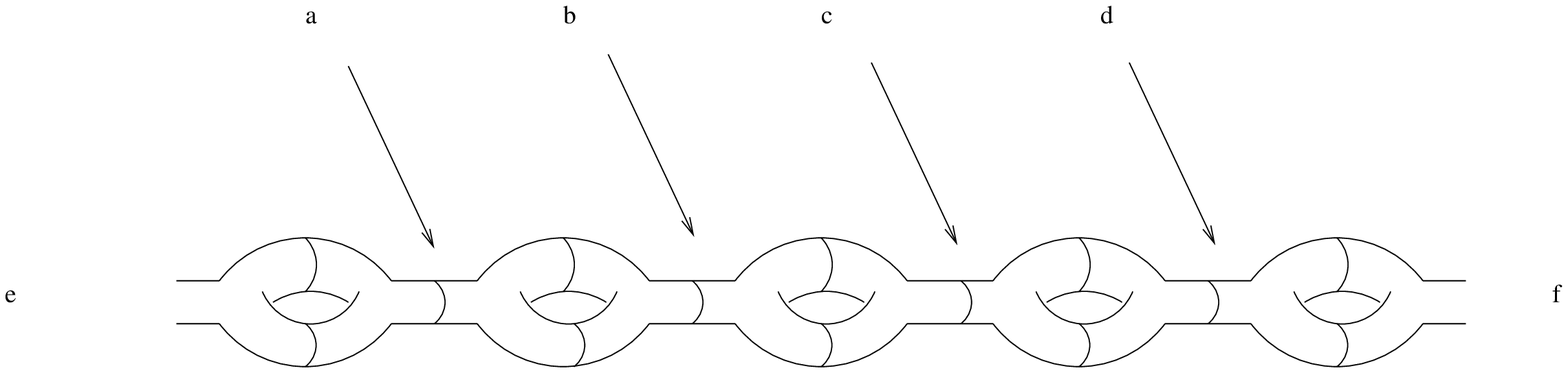}
\caption{\small {}}
\label{bi-infinite}
\end{figure}
\bigskip

\begin{example}\label{ex:ls1} As a different example of the same phenomenon, we take  $X$ to be the surface represented in Figure \ref{bi-infinite} and we let $T$ be the step-1 right-translation suggested in this figure, transforming $C_i$ into $C_{i+1}$ for every $i$ in $\mathbb{Z}$. Consider a hyperbolic structure $H_0$ on $X$ that is translation-invariant. In other words, we suppose that all the two-holed tori that are in the complement of the curves $C_i$ are isometric, and that the twist parameters about these  curves are all equal. For each $i=1,2,\ldots$, let $\tau_i$ be the $i$-th power of the positive Dehn twist along $C_i$. Let $f=\tau_1\circ\tau_2\circ\tau_3\circ\ldots$ be  the infinite composition of the $\tau_i$'s and let
$H_1$ be the hyperbolic structure $f(H_0)$. By using the Hopf-Rinow theorem as in Example \ref{ex:ls} above, we can see that the metrics $H_0$ and $H_1$ are complete. Let $D_0$ be an essential simple closed curve on $X$ satisfying $i(C_0,D_0)>0$
and $i(C_i,D_0)=0$ for all $i\not=0$. For each $i$ in $\mathbb{Z}$, let $D_i=T^i(D_0)$. Then, we have \[\frac{l_{H_{1}}(D_i)}{l_{H_{0}}(D_i)}\to\infty\] as $i\to\infty$. Therefore, $H_0$ and $H_1$ are not related by a length-spectrum bounded homeomorphism. Hence, $H_1\not\in\mathcal{T}_{ls}(H_0)$, and again we have $\mathcal{T}_{ls}(H_0)\not=\mathcal{T}_{ls}(H_1)$.
\end{example}

More generally, we have the following:
\begin{proposition}\label{prop:eq}
Given a hyperbolic structure $H_0$ on an infinite type surface $X$, there exists a hyperbolic structure $H_1$ on $X$ such that $H_0$ and $H_1$ are not length-spectrum equivalent. In particular, $\mathcal{T}_{ls}(H_0)\not=\mathcal{T}_{ls}(H_1)$.
\end{proposition}
\begin{proof}
Consider an infinite collection $C_i$, $i=1,2,\ldots$ of disjoint pairwise non-homotopic essential curves on $X$.  Up to taking a subsequence of $(C_i)$, we can find, for each  $i=1,2,\ldots$, a curve $D_i$ satisfying $i(C_i,D_i)>0$
and $i(C_i,D_j)=0$ for all $j\not=i$.
For each  $i=1,2,\ldots$, we take $\tau_i:X\to X$ to be a power of a positive Dehn twist along $C_i$ satisfying $l_{H_{0}}(\tau_i(D_i))>i$. Let $f$ be the infinite composition $f=\tau_1\circ\tau_2\circ\ldots$. Then,
 \[\frac{l_{H_{0}}(f(D_i))}{l_{H_{0}}(D_i)}\to\infty\] as $i\to\infty$.
This shows that the hyperbolic structure $f(H_0)$ is not related to $H_0$ by a length-spectrum bounded homeomorphism.
\end{proof}

            We denote by $\mathrm{Hom}(X)$ the group of self-homeomorphisms of $X$, and by
        $\mathrm{Hom}_{ls}(H_0)$ the subgroup of $\mathrm{Hom}(X)$ consisting of self-homeomorphisms that are length-spectrum bounded with respect to the metric $H_0$. We observe that in general $\mathrm{Hom}_{ls}(H_0)\not=\mathrm{Hom}(X)$. This can be seen in the following example.

\bigskip
  \begin{figure}[!hbp]
\centering
\psfrag{a}{\small $C_1$}
\psfrag{b}{\small $C_2$}
\psfrag{c}{\small $C'_1$}
\psfrag{d}{\small $C'_2$}
\psfrag{e}{\small $C''_1$}
\psfrag{f}{\small $C''_2$}
\psfrag{g}{\small $\ldots$}
\includegraphics[width=.65\linewidth]{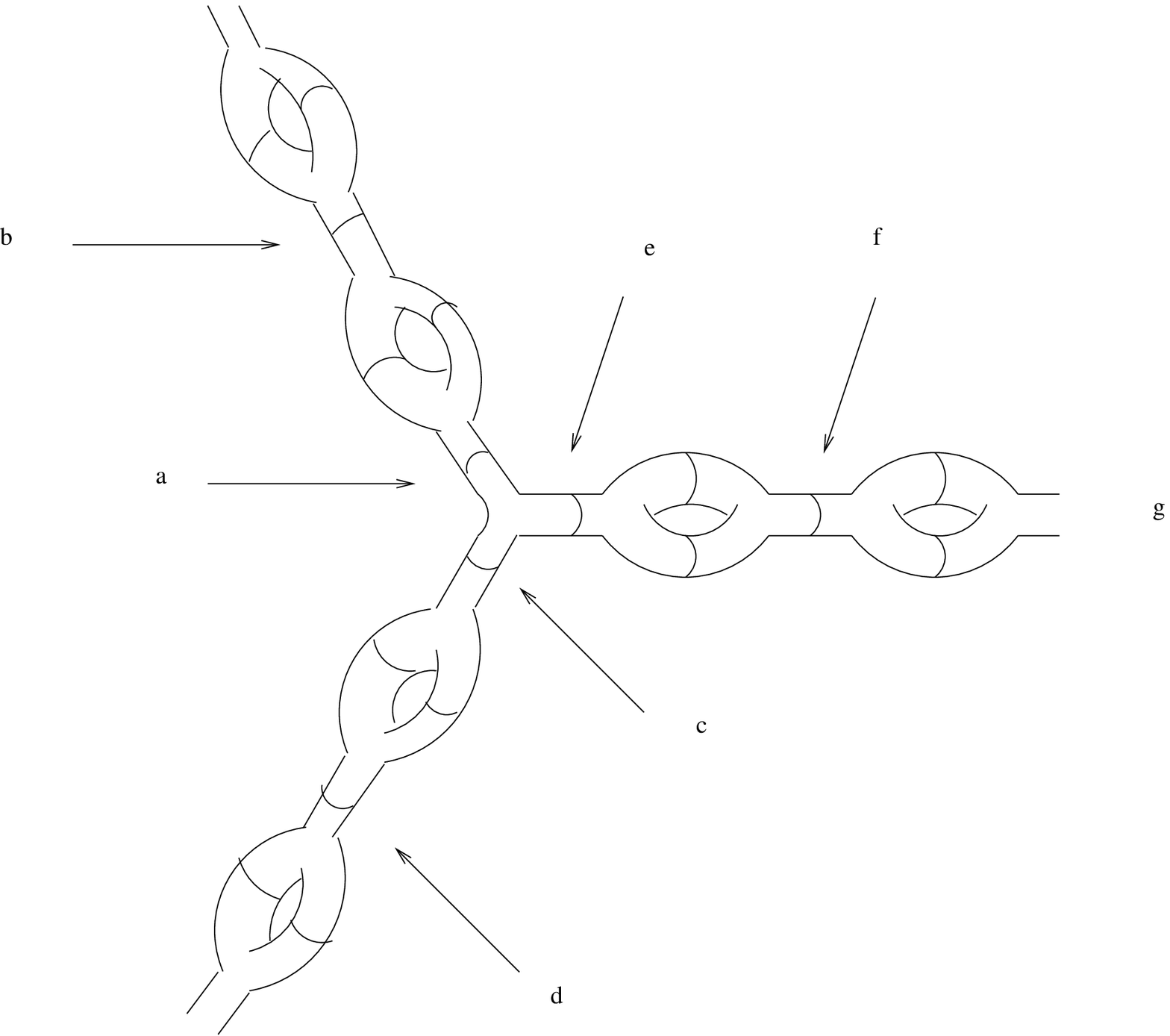}
\caption{\small {}}
\label{tri-infinite}
\end{figure}
\bigskip

        \begin{example}\label{ex:tripod} Let $X$ be the topological surface represented in Figure \ref{tri-infinite}. We consider a homeomorphism $f$ obtained by applying a $2\pi/3$-counterclockwise rotation suggested by the picture. For $i=1,2,\ldots$, the curves $C_i$, $C'_i$ and $C''_i$  are in the homotopy classes represented in this figure. We choose a hyperbolic metric $H_0$ on $X$  such that for all $i=1,2,\ldots$,  $l_{H_{0}}(C_i)=1/i$, $l_{H_{0}}(C'_i)=i$ and $l_{H_{0}}(C''_i)=i^2$. By using, as in the previous examples, the theorem of Hopf-Rinow, we can insure that this metric is complete.
 It is clear that $f$ is not homotopic to a length-spectrum bounded homeomorphism.  
        \end{example}

        More generally, we have the following.
        \begin{proposition}\label{prop:ne} For any hyperbolic structure $H_0$ on a surface $X$ of infinite type, we have  $\mathrm{Hom}_{ls}(H_0)\varsubsetneq\mathrm{Hom}(X)$.
        \end{proposition}

        \begin{proof} The homeomorphism of $f:X\to X$ constructed in the proof of Proposition \ref{prop:eq} is not in  $\mathrm{Hom}_{ls}(H_0)$.
        \end{proof}

         We let $\mathrm{MCG}(X)$ denote the quotient of the group $\mathrm{Hom}(X)$ by the normal subgroup consisting of orientation-preserving homemorphisms that are homotopic to the identity. (Recall that our homotopies preserve setwise the boundary components, at each time.) The group $\mathrm{MCG}(X)$ is the {\it mapping class group} of $X$.

    \begin{definition}[The  length-spectrum mapping class group] The {\it length-spectrum mapping class group of $X$ relative to the hyperbolic metric $H_0$}, which we denote by $\mathrm{MCG}_{ls}(H_0)$, is the subgroup of $\mathrm{MCG}(X)$ consisting of homotopy classes of homeomorphisms that are length-spectrum bounded with respect to the metric $H_0$.
        \end{definition}

      It follows from Proposition \ref{prop:ne} that for any hyperbolic structure $H_0$ on a surface $X$ of infinite type, we have  $\mathrm{MCG}_{ls}(H_0)\varsubsetneq\mathrm{MCG}(X)$.

        The group $\mathrm{MCG}_{ls}(H_0)$ acts on the Teichm\"uller space $\mathcal{T}_{ls}(H_0)$ by the quotient of the action
of homeomorphisms on marked hyperbolic surfaces by right-composition by the inverse. More precisely, for each marked hyperbolic surface $(f,H)$ and for each homeomorphism $g:(X,H_0)\to (X,H)$ which is  length-spectrum bounded, we set
\begin{equation}\label{eq:action-pre}
g\cdot (f,H)= (f\circ g^{-1}, H).
\end{equation}
     This action preserves the metric $d_{ls}$ on $\mathcal{T}_{ls}(H_0)$.

The following is a natural question:
\begin{question}\label{q:im}
Given any hyperbolic metric on an infinite-type surface, is  the natural map from  $\mathrm{MCG}_{ls}(H_0)$ into the isometry group $\mathrm{Isom}(\mathcal{T}_{ls}(H_0),d_{ls})$ of the metric space $(\mathcal{T}_{ls}(H_0),d_{ls})$ an isomorphism ? Is it an isomorphism onto an index-2 subgroup ? (The reason for talking about an index-two subgroup is that one may have to extend the image of $\mathrm{MCG}_{ls}(H_0)$ in $\mathrm{Isom}(\mathcal{T}_{ls}(H_0),d_{ls})$ by the homotopy class of an orientation-reversing length-spectrum bounded homeomorphism of $(X,H_0)$, if such a homeomorphism exists.)
\end{question}

\begin{remark}
The study of the image of  $\mathrm{MCG}_{ls}(H_0)$ in the isometry group $\mathrm{Isom}(\mathcal{T}_{ls}(H_0),d_{ls})$ can also be asked in the case of finite-type surfaces, and it does not seem to have been settled. We recall that in the case of a finite-type hyperbolic surface $S$, we have $\mathrm{MCG}_{ls}(S)\simeq \mathrm{MCG}_{}(S)$, and that
a famous result by Royden, completed by Earle and Kra, says that the natural homomorphism from the extended mapping-class group of $S$ (which is an order-two extension of  $\mathrm{MCG}_{}(S)$) into the isometry group of $\mathcal{T}_{}(S)$ equipped with the Teichm\"uller metric is an isomorphism, except for some finite set of surfaces of small genus and small number of boundary components (see \cite{Ro} and \cite{EK}).
 For surfaces of finite type, Li and Liu (see \cite {Li1986} and  \cite{Liu1999}) proved  that the length-spectrum metric on Teichm\"uller space defines the same topology as the Teichm\"uller metric, but of course this result is not sufficient to prove that the result of Royden and Earle-Kra also applies to the length-spectrum metric.
\end{remark}

            We denote by  $\mathrm{MCG}_{f}(X)$ the subgroup of $\mathrm{MCG}(X)$ consisting of the homotopy classes of homeomorphisms supported on a subsurface of $X$ of finite topological type.

            \begin{proposition}\label{prop:inclusion}
            We have a natural inclusion  $\mathrm{MCG}_{f}(X)\subset  \mathrm{MCG}_{ls}(H_0)$.
            \end{proposition}

\begin{proof}

     Since every element of the mapping class group of a surface of finite topological type is a finite composition of Dehn twists, it suffices to show that a Dehn twist about a simple closed curve in $X$ is
     length spectrum bounded.

     Let $\tau_{\beta}$ be a Dehn twist about a simple closed curve
     $\beta$ and let $\alpha$ be an arbitrary element of $\mathcal{S}(X)$. 
     
      If $i(\alpha, \beta)=0$, then
     $\alpha=\tau_{\beta}(\alpha)$ and
     $l_{H_0}(\alpha)=l_{H_0}(\tau_{\beta}(\alpha))$.

     If  $i(\alpha,\beta) \ne 0$, then, by the Collar Lemma (cf.
 \cite{Buser}), there exists a
     constant $t_0$, depending only on the hyperbolic length
     $l_{H_0}(\beta)$ such that the cylinder in $X$ defined as a $t_0/2$-neighborhood of the closed geodesic $\beta$ is embedded in $X$.
    Since $\alpha$ traverses this cylinder $i(\alpha,\beta)$ times, we hace, $l_{H_0}(\alpha) \ge i(\alpha, \beta) t_0$. Furthermore, from the definition of a Dehn twist, we have
     \[l_{H_0}(\tau_{\beta}(\alpha)) \le l_{H_0}(\alpha)+i(\alpha,
     \beta) l_{H_0}(\beta).\] This gives
     $$\frac{l_{H_0}(\tau_{\beta}(\alpha))} {l_{H_0}(\alpha)} \le 1+
     \frac{l_{H_0}(\beta)} {t_0}.
     $$

     Similarly, we have
    $$\frac{l_{H_0}(\alpha)} {l_{H_0}(\tau_{\beta}(\alpha))}\le 1+
     \frac{l_{H_0}(\beta)} {t_0}.
     $$

     This implies that
\[\log \sup_{\alpha\in\mathcal{S}(X)} \left\{
            \frac{l_{H_0}(\tau_{\beta}(\alpha))}{l_{\phi(H_0)}(\alpha)},
            \frac{l_{\phi(H_0)}(\alpha)}{l_{H_0}(\tau_{\beta}(\alpha))}\right\}<\infty.\]

Thus, $\tau_{\beta}$ is length spectrum bounded.

     \end{proof}

Note that the inclusion in Proposition \ref{prop:inclusion} is in general strict.

\begin{proposition}\label{prop:t0}
Let $H$ be a hyperbolic metric on $X$ such that there exists a
sequence $\alpha_n$, $n=1,2,\ldots$ of elements  in
$\mathcal{S}(X)$ represented by disjoint simple closed curves
satisfying $l_H(\alpha_n)\to 0$ as $n\to\infty$. For every
$n=1,2,\ldots$, let $\tau_n$ be the positive Dehn twist about
$\alpha_n$. Then, $$d_{ls}(\tau_n(H),H)\to 0.$$
\end{proposition}

\begin{proof}
Note that $l_{\tau_{n}(H)}(\alpha)= l_H(\tau_n^{-1}(\alpha))$.
  We need to prove that
\begin{equation}\label{eq:tends}\log \sup_{\alpha\in\mathcal{S}(X)} \left\{
\frac{l_{H}( \tau_n^{-1}(\alpha))}{l_{H}(\alpha)},
\frac{l_{H}(\alpha)}{l_{H}(\tau_{n}^{-1}(\alpha))}\right\}\to 0 \hbox{ as } n\to\infty.
\end{equation}
Let $\alpha$
be an element of $\mathcal{S}(X)$ and let  $n$ be a positive integer.

If $i(\alpha, \alpha_n)=0$,
  then $\alpha= \tau_n^{-1}(\alpha)$ and  $l_{H} (\tau_n^{-1}(\alpha))= l_{H}
(\alpha)$.

Assume that $i(\alpha, \alpha_n) \neq 0$ and let $\varepsilon_n=l_H(\alpha_n)$.
By a version of
the Collar Lemma \cite{Buser}, there exists $B_{\alpha}>0$ such that
\[l_{H}(\alpha)=   i(\alpha, \alpha_n) |\log\varepsilon_n|+ B_{\alpha}.\]

From the definition of a Dehn twist, we have
\[ i(\alpha, \alpha_n) |\log\varepsilon_n| + B_{\alpha} -
 i(\alpha, \alpha_n) \varepsilon_n \leq
l_{H}(\tau_n^{-1}(\alpha))\]
\[\leq  i(\alpha, \alpha_n)
| \log\varepsilon_n| + B_{\alpha} +  i(\alpha, \alpha_n)
\varepsilon_n .\]

  Then we have
  \begin{eqnarray*}
\frac{l_{H}(\tau_n^{-1}(\alpha))}{l_{H}(\alpha)}&\leq&
 \frac{ i(\alpha, \alpha_n) |\log\varepsilon_n| +
B_{\alpha} + i(\alpha, \alpha_n) \varepsilon_n} {  i(\alpha,
\alpha_n) |\log\varepsilon_n|+ B_{\alpha}}
 \\&=&
  1+ \frac{ i(\alpha, \alpha_n) \varepsilon_n}{ i(\alpha, \alpha_n) |\log\varepsilon_n|
  + B_{\alpha}}
  \\&\leq&
   1+ \frac{ i(\alpha, \alpha_n)\varepsilon_n}{ i(\alpha, \alpha_n) |\log
\varepsilon_n|}
  \\&\leq&
 1+   \frac {\varepsilon_n} {|\log
\varepsilon_n|} .\end{eqnarray*}

In the same way, we can show that
\[\frac{l_{H}(\alpha)}{l_{H}(\tau_n^{-1}(\alpha))}\leq  1+   \frac {\varepsilon_n} {|\log
\varepsilon_n|} .\]

This proves that

\[ \sup_{\alpha\in\mathcal{S}(X)} \left\{
\frac{l_{H}(\tau_n^{-1}(\alpha))}{l_{H}(\alpha)},
\frac{l_{H}(\alpha)}{l_{H}(\tau_{n}^{-1}(\alpha))}\right\}\to 1 \hbox{ as } n\to\infty,\]
which implies (\ref{eq:tends}).
\end{proof}

We shall say that the action of a group on a topological space is {\it discrete} if the orbit of any point is discrete.

As a consequence of Proposition \ref{prop:t0}, we have the following.
\begin{corollary} There exist hyperbolic structures $H$ such that the action of the group $\mathrm{MCG}_{ls}(H)$ on the Teichm\"uller
space $\mathcal{T}_{ls}(H)$ equipped with the metric $d_{ls}$ is not discrete \end{corollary}

\begin{corollary}
A necessary condition for the action of the group
$\mathrm{MCG}_{ls}(H)$ on the Teichm\"uller space
$\mathcal{T}_{ls}(H)$ equipped with the metric $d_{ls}$ to be
discrete is that there exists a 
constant $\epsilon>0$, such that for any element $\alpha$ in
$\mathcal{S}(X)$, $l_{H}(\alpha)\ge \epsilon$.
\end{corollary}

\begin{question}
Find necessary and sufficient conditions on a hyperbolic structure $H$ on a surface of infinite type such that the action of  the group $\mathrm{MCG}_{ls}(H)$ on the Teichm\"uller space
$\mathcal{T}_{ls}(H)$ equipped with the metric $d_{ls}$ is discrete. (It is known that for surfaces of finite analytical type, the action is always  discrete.)
\end{question}

\begin{remark}
It is known that for surfaces of infinite analytical type, the action of the quasiconformal mapping class group $\mathrm{MCG}_{qc}(H)$ on the quasiconformal Teichm\"uller space $\mathcal{T}_{qc}(H_0)$ is in general not discrete (see Section \ref{s11} below for the definitions, and \cite{M2003}, \cite{F2004a} and \cite{Fujikawa} for work on this question).
\end{remark}

    Proposition \ref{prop:eq} implies that for every hyperbolic metric $H_0$ on an infinite-type surface, we have
            $\mathrm{MCG}_{ls}(H_0)\not=\mathrm{MCG}(X)$.

\begin{example}[A noncompact closed ball] \label{ex:proper} Let $X$ be the topological surface represented in Figure \ref{bi-infinite}. We let $f:X\to X$ be the homeomorphism defined as the step-one translation to the right suggested by this picture, and we equip $X$ with a hyperbolic metric $H$ which is invariant by the homeomorphism $f$. In particular, the homotopy classes of curves $C_n$ satisy  $C_n=f^{n}(C_0)$ for every $n\in\mathbb{Z}$.

Let $x\in\mathcal{T}_{ls}(H)$ denote the equivalence class of the marked hyperbolic surface $(\mathrm{Id}, H)$. For each $n\in\mathbb{Z}$, let $\tau_n$ denote the positive Dehn twist along $C_n$ and let $x_n\in\mathcal{T}_{ls}(H)$ denote the equivalence class of the marked surface $(\tau_n,H)$.

From the translation-invariance of $H$, the distance $d_{ls}(x,x_n)$
does not depend on $n$. Thus, all the points $x_n$ are contained in
a geometric ball $B$ centered at $x$ and of a fixed radius. 

We claim that there
exists a constant $K>0$ such that
\begin{equation} \label{eq:K} \forall m\not= n\in\mathbb{Z},\ d_{ls}(x_m,x_n)\geq K.
\end{equation}
Indeed, let $\alpha_0$ be the homotopy class of an essential curve satisfying $i(\alpha_0,C_0)\not=0$ and $i(\alpha_0,C_n)=0$ for all $n\not=0$ and for all $m$ in $\mathbb{Z}$, let $\alpha_m$ denote the homotopy class $f^m(\alpha_0)$. Then, for all integers $m\not=n$, the value 
\[\displaystyle\max \left\{\log \frac{l_{H_{m}}(\alpha_m)}{l_{H_{n}}(\alpha_m)}, \log \frac{l_{H_{n}}(\alpha_m)}{l_{H_{m}}(\alpha_m)}\right\}\] is independent of $m$ and $n$. Calling this constant $K$, we have (\ref{eq:K}).

Property (\ref{eq:K}) implies that the sequence $x_n$ does not have a convergent subsequence. Therefore the closed ball $B$ is not compact.

\end{example}

We recall that a metric space is called {\it proper} if every closed ball is compact.
 From the preceding example, there exist hyperbolic structures $H$ such that the Teichm\"uller space  $\mathcal{T}_{ls}(H)$ is not proper. We shall study properness in more generality in Proposition \ref{prop:notproper} below. For that, we shall use the following.

        \bigskip
  \begin{figure}[!hbp]
\centering
\psfrag{a}{\small $\alpha$}
\psfrag{b}{\small $\beta$}
\includegraphics[width=.35\linewidth]{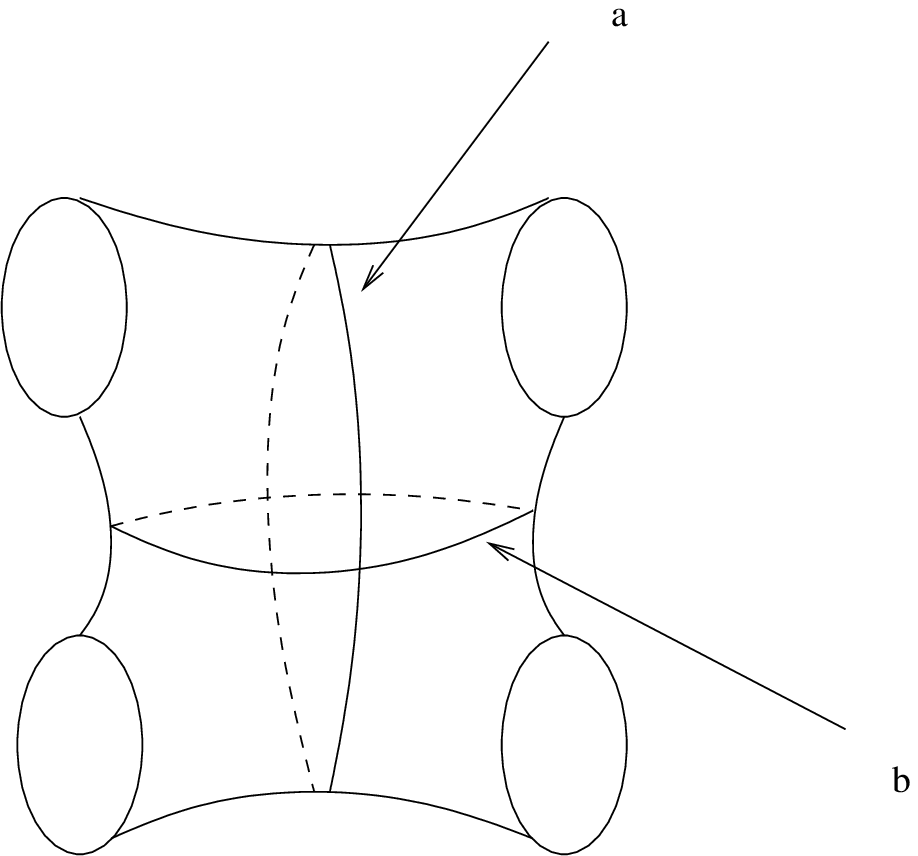}
\caption{\small {}}
\label{fig:sphere}
\end{figure}
\bigskip

 Let $F$ be a compact sphere with four holes, equipped with a hyperbolic structure $H$ with geodesic boundary and suppose that $F$ is decomposed into a union of two hyperbolic pairs of pants, $P_1$ and $P_2$. (See Figure \ref{fig:sphere} in which $\beta$ is a boundary components of $P_1$ and of $P_2$.) We assume that there exists a positive constant $M$ such that each boundary component $C$ of either of the two pairs of pants $P_1$ and $P_2$, including the curve $\beta$ itself, satisfies
 \[\frac{1}{M}\leq l_{H}(C)\leq M.\]
We let $\alpha$ be a simple closed geodesic on $F$ defined as follows. We first take, in each of the pairs of pants $P_1$ and $P_2$ the shortest geodesic segment  joining $\beta$ to itself and not homotopic to a point relative to the curve $\beta$ (this shortest geodesic segment is called the {\it seam} of the pair of pants). Then, we join together the endpoints of these two seams by using geodesic segments contained in $\beta$, producing a closed curve in
$F$. Then, the curve $\alpha$ is the simple closed geodesic homotopic to that curve.
 (There are two possibilities for joining the endpoints of the two seams, namely by turning ``to the left" or ``to the right" along $\beta$, and making such a choice will not affect the result we prove in the next lemma.)

\begin{lemma} \label{lem:four}
With the above notation, there exists a constant $K>1$ that depends
only on $M$ such that if $\tau$ denotes the positive Dehn
twist along $\beta$, we have
\[\max \left(\frac{l_{H}(\tau^{-1}(\alpha))}{l_{H}(\alpha)},
\frac{l_{H}(\alpha)}{l_{H}(\tau^{-1}(\alpha))},
\frac{l_{H}(\tau(\alpha))}{l_{H}(\alpha)},
\frac{l_{H}(\alpha)}{l_{H}(\tau(\alpha))}\right) >K.\]

 \end{lemma}

  \begin{proof}
 The isometry type of the sphere with four holes $F$ is determined by the values of six parameters, five of them being the lengths of the geodesic boundary components of the pairs of pants $P_1$ and $P_2$, of which two are equal, and the sixth value being the  twist parameter along the curve $\beta$. The values $l_H(\tau(\alpha))$ and  $l_H(\tau^{-1}(\alpha))$ are continuous functions of these parameters. Furthermore, for any element $\alpha$ in $\mathcal{S}(F)$, at least one of the two values $\vert l_H(\tau(\alpha))-l_H(\alpha)\vert $ and
 $\vert l_H(\tau^{-1}(\alpha))-l_H(\alpha)\vert $ is nonzero. This can be deduced from the fact that the length function $\theta\mapsto l_H(\tau_{\theta}(\alpha))$ under the twist deformation along the curve $\beta$ is strictly convex (see \cite{FLP} p. 130).

 Therefore, for each $\alpha\in\mathcal{S}$, the value
\[ \max \left(\frac{l_{H}(\tau^{-1}(\alpha))}{l_{H}(\alpha)},
\frac{l_{H}(\alpha)}{l_{H}(\tau^{-1}(\alpha))},
\frac{l_{H}(\tau(\alpha))}{l_{H}(\alpha)},
\frac{l_{H}(\alpha)}{l_{H}(\tau(\alpha))}\right) \] is strictly positive. Since in this maximum, each time a number occurs, then its inverse also occurs, the value of the maximum is $>1$. This value varies continuously in terms of the six parameters referred to above. Under the hypotheses considered, the six parameters vary in a compact set. This shows the existence of the constant $K$.
 \end{proof}

We now consider a surface $X$, equipped with a hyperbolic structure $H$ and a pair of pants decomposition $\mathcal{P}$ satisfying Property ($\star$) stated in the introduction, and such that if
 $C_i$, $i=1,2,\ldots$ is the collection of simple closed curves that are boundaries of pairs of pants in the collection $\mathcal{P}$, then,
  \begin{equation}\label{LM} \exists M>0,\forall i=1,2,\ldots, \frac{1}{M}\leq l_H(C_i)\leq M.
 \end{equation}

\begin{proposition}\label{prop:notproper} Let $H$ be a hyperbolic metric on an infinite type surface $X$ satisfying Property (\ref{LM}). Then, the Teichm\"uller space  $\mathcal{T}_{ls}(H)$ is not proper.
\end{proposition}

\begin{proof}

We can find a collection $F_i$, $i=1,2,\ldots$ of four-holed spheres with geodesic boundary embedded in $X$ and satisfying the following.
\begin{enumerate}
\item for each $i=1,2,\ldots$, the four-holed sphere $F_i$ is the union of two pairs of pants in the collection $\mathcal{P}$;
\item for all $i\not=j$, the intersection of the interiors of $F_i$ and $F_j$ is empty.
\end{enumerate}
In particular, in the interior of each four-holed sphere $F_i$, there is a curve $\beta$ which belongs to the collection $\{C_i\}$. For each four-holed sphere $F_i$, we denote by $\tau_i$ the positive Dehn twist along the curve $\beta$ in its interior. We let $x_i\in\mathcal{T}_{ls}(H)$ be the point representing $(\tau_i,H)$ and we let $x$ be the point representing $(\mathrm{Id},H)$. By the arguments used in the proof of Proposition \ref{prop:inclusion}, there exists a constant $R>0$ such that we have
\[\forall i=1,2,\ldots,\ d_{ls}(x,x_i)\leq R.\]
Thus, the sequence $(x_i)$ is contained in the closed ball $B$ of radius R and centered at $x$.
From Lemma \ref{lem:four}, we deduce that there exists a constant $K$ such that for each $i\not=j$, $d_{ls}(x_i,x_j)\geq K$. This shows that the sequence $x_i$ does not have any Cauchy subsequence. This implies that the closed ball $B$ is not compact, which proves that $\mathcal{T}_{ls}(H)$ is not proper.
\end{proof}

It should be true that for any hyperbolic structure $H$ on an infinite type surface, $\mathcal{T}_{ls}(H)$ is not proper, and this should follow from the fact that this space is infinite-dimensional.

We also ask the following.

\begin{question}\label{q:3}
If $H$ and $H'$ are arbitrary hyperbolic metrics on surfaces on an
infinite-type surface $X$, are the two Teichm\"uller spaces
$\mathcal{T}_{ls}(H)$ and   $\mathcal{T}_{ls}(H')$ locally
bi-Lipschitz equivalent ?
\end{question}

In this regard, we note that in \cite{Fletcher}, A. Fletcher proved that the answer to Question \ref{q:3} is yes if instead of the  length-spectrum metric, one takes the Teichm\"uller metric (see Section \ref{s11} below). As a matter of fact, Fletcher proved that the Teichm\"uller metric on any Teichm\"uller space of an infinite-type Riemann surface is locally bi-Lipschitz equivalent to the Banach space $l^{\infty}$ of bounded sequences  with the supremum norm.

We end this section by a discussion of the completeness of the length spectrum metric.

 In the paper \cite{Shiga2003}  (Corollary 1.1), H. Shiga showed that for surfaces of infinite type, the topology defined by the length spectrum metric is in general not complete, but that under condition (\ref{LM}), this metric is complete. A natural question then is the following:
 \begin{question} Give a necessary and sufficient condition on a hyperbolic metric $H$ on a surface of  infinite type under which the Teichm\"uller space  $\mathcal{T}_{ls}(H)$ is complete.
\end{question}

Note that for
surfaces of finite type, the completeness of the length spectrum
metric does not follow automatically from the completeness of the Teichm\"uller
metric, even though the topologies induced by the length spectrum and
the Teichm\"uller metrics are the same (see \cite{Liu1999}). But it is true that the length spectrum metric in the case of surfaces of finite type is complete. For completeness, we now provide a proof of this fact.

We shall discuss the Teichm\"uller metric for Teichm\"uller spaces of surfaces of infinite type in Section \ref{s11} below, but we assume here that the reader is familiar with the Teichm\"uller metric of surfaces of finite type.

 Let $S$ be a  Riemann surface of conformally finite type $(g,n)$, where $g$ is the genus and $n$ the
number of punctures of $S$. For any $\epsilon >0$, let
$\mathcal{T}^{\epsilon}(S)$ be the subset of $\mathcal{T}(S)$  consisting of
equivalence classes of marked hyperbolic surfaces in which the length of any simple closed curve which
is not homotopic to a puncture is at least $\epsilon$. We call
$\mathcal{T}^{\epsilon}(S)$ the {\it thick part} of  $\mathcal{T}(S)$. 
We recall that for surfaces of finite type, all the Teichm\"uller spaces coincide from the set-theoretic point of view, and therefore we do not distinguish between  the spaces $\mathcal{T}_{ls}$ and $\mathcal{T}_{qc}$. We denote by 
$d_{ls}$ and $d_{qc}$ respectively the length-spectrum and the Teichm\"uller metrics on $\mathcal{T}(S)$.

We shall use the following result, proved in \cite{Liu2008}.
\begin{lemma} \label{lem:Liu}
 There  exists  a constant $N(g,n,
\epsilon)$, depending only on $g, n, \epsilon$, such that for any $x$ and $y$ in
$\mathcal{T}^{\epsilon}(S)$, we have
 $$
 d_{ls}(x,y) \le d_{qc}(x,y)
 \le 4 d_{ls}(x,y)+N(g,n, \epsilon).
 $$
 \end{lemma}

We have the following.

\begin{theorem}\label{th:LL}  Let $S$ be a conformally finite type hyperbolic surface and let $x_n, n=0, 1, \ldots$ be a sequence of elements in $\mathcal{T}(S)$.
Then 
\[\lim_{n \to \infty}d_{qc}(x_n, x_0)= \infty\iff\lim_{n \to \infty}d_{ls}(x_n, x_0)= \infty.\]
 \end{theorem}
 \begin{proof} We first recall that a result of Wolpert (see \cite{Abikoff}) says that 
  given any $K$-quasiconformal map $f:S_0\to S$, then, for any hyperbolic metrics $H_0$ and $H$ on $S_0$ and $S$ respectively and for any element $\alpha$ in $\mathcal{S}(S_0)$, we have
\begin{equation}\label{eq:Wolpert}
\frac{1}{K}\leq\frac{l_H(f(\alpha))}{ l_{H_{0}}(\alpha)}\leq K.
\end{equation} 
This implies that if $\lim_{n \to
 \infty}d_{ls}(x_n, x_0)=\infty$, then
  $\lim_{n \to
 \infty}d_{qc}(x_n, x_0)=\infty$.

 Now suppose that $\lim_{n \to  \infty}d_{qc}(x_n, x_0)=\infty$. Consider first the following two cases:
 
  (1) the sequence $(x_n)$ leaves any thick part  $\mathcal{T}^{\epsilon}(S)$. In
 this case it is clear that $\lim_{n \to
 \infty}d_{ls}(x_n, x_0)=\infty$. 
 
 (2) the sequence $(x_n)$ stays in some thick part
 $\mathcal{T}^{\epsilon}(S)$ for some $\epsilon>0$. In
 this case, from Lemma \ref{lem:Liu} we have $\lim_{n \to
 \infty}d_{ls}(x_n, x_0)=\infty$.
 
 Now we discuss the  general case.  Assume that $d_{ls}(x_n, x_0)$ does not go to $ \infty$ as $n\to \infty$. Then there exists subsequence
 $x_{n_k}$, $n_k=1,2,, \ldots$ of $(x_n)$ such that the set $\{d_{ls}(x_{n_k},
 x_0), n_k=1, \cdots\}$ is bounded. From the above discussion, if for some $K_0$, a subsequence $x_{n_k}$,
 $ n_k>K_0$ leaves any thick part of
 $\mathcal{T}^{\epsilon}(S)$, we get a contradiction.
 Otherwise, there exists a subsequence of $x_{n_k}$,$ n_k=1,2,
 \ldots$  which stays in some thick part of
 $\mathcal{T}^{\epsilon}(S)$, which also leads a contradiction.
 \end{proof}

 \begin{theorem}\label{th:co} Let $S$ be a conformally finite type Riemann surface. Then $\mathcal{T}_{ls}(S)$ equipped with the length spectrum metric is complete.
 \end{theorem}

 \begin{proof} We need to
 prove that for any Cauchy sequence $x_n, n=1, \ldots$ in $(\mathcal{T}(S),d_{ls})$, there exists $x_0 \in \mathcal{T}(S)$, such that $\lim_{n \to \infty}
 d_{ls}(x_n, x_0)=0$.

The sequence $(x_n)$ is bounded in the metric $d_{ls}$. From Theorem \ref{th:LL}, it follows that $(x_n)$ is also bounded in the metric $d_{qc}$. Since $d_{qc}$ is complete and since the Teichm\"uller space considered is finite-dimensional, there exists a subsequence $(x_{n_k})$ of $(x_n)$ and $x_0 \in \mathcal{T}(S)$ such that
 $\lim_{n_k \to \infty} d_{qc}(x_{n_k}, x_0)=0$.  From the topological
 equivalence of $d_{qc}$ and $d_{ls}$ on $\mathcal{T}(S)$
 \cite{Liu1999}, it follows that
 $$\lim_{n_k \to \infty} d_{ls}(x_{n_k},
 x_0)=0.$$
 As $\{x_n
 \}^{\infty}_{n=1}$ is a Cauchy sequence in the metric $d_{ls}$,
 $\lim_{n \to \infty} d_{ls}(x_{n},
 x_0)=0$.
 \end{proof}

From \cite{Shiga2003} and \cite{Liu2008},  Theorem \ref{th:co}  is not true for certain Teichm\"uller spaces of infinite type
 surfaces.

        \section{The bi-Lipschitz Teichm\"uller spaces} \label{s22}

        We say that a homeomorphism $f:(X_1,d_1)\to (X_2,d_2)$ between metric spaces is  {\it bi-Lipschitz} if there exists a real number $K\geq 1$ satisfying
        \begin{equation}\label{eq:bl}\frac{1}{K} d_1(x,y)\leq d_2 (f(x),f(y))\leq Kd_1(x,y).
        \end{equation}
        The real number $K$ in (\ref{eq:bl}) is called a {\it bi-Lipschitz constant of $f$}.
        Two metric spaces are said to be {\it bi-Lipschitz equivalent} if there exists a bi-Lipschitz homeomorphism between them.

        \begin{definition}  Let $H_0$ be a hyperbolic metric on the surface $X$. The {\it bi-Lipschitz Teichm\"uller space} of $H_0$, denoted by $\mathcal{T}_{bL}(H_0)$, is the set of equivalence classes $[f,H]$ of pairs $(f,H)$ where $H$ is a hyperbolic metric on a surface homeomorphic to $X$ and where $f:H_0\to H$ (the {\it marking} of $H$)  is a bi-Lipschitz homeomorphism, and where two such pairs  $(f,H)$ and $(f',H')$ are considered equivalent if there exists an isometry $f'':H\to H'$ homotopic to $f'\circ  f^{-1}$.
         \end{definition}

         The topology of $\mathcal{T}_{bL}(H_0)$ is the one induced by the {\it bi-Lipschitz} metric $d_{bL}$, defined by 
                 \[d_{bL}([f,H],[f',H'])=\frac{1}{2}\log \inf \{K\},\]
                  where the infimum is taken over
        all bi-Lipschitz constants $K$ of
        bi-Lipschitz  homeomorphisms $f'':H\to H'$  homotopic to $f'\circ  f^{-1}$, with $(f,H)$ and $(f',H')$ being two hyperbolic surfaces representing the two points
         $[f,H]$ and $[f',H']$   in
 $\mathcal{T}_{bL}(H_0)$.
 
  To see that this defines indeed a metric on
$\mathcal{T}_{bL}(H_0)$., we note that if we take a sequence of bi-Lipschitz maps homotopic to $f'\circ f^{-1}$ whose bi-Lipschitz constants tend to 1, then there exists an isometry homotopic to  $f'\circ f^{-1}$. One easy way to see this is to use the following remarks that show that under these hypotheses, the length-spectrum constant for the map  $f'\circ f^{-1}$ must be 1, hence, by Proposition \ref{prop:isometry-}, there exists an isometry homotopic to  $f'\circ f^{-1}$.

        If $f:H_0\to H$ is a bi-Lipschitz homeomorphism with constant $K$, then, from the definition of the length of a curve, we have, for every element $\alpha$ in $\mathcal{S}_0$,
        \[\frac{1}{K}\leq \frac{l_H(f(\alpha))}{l_H(\alpha)}\leq K.\]
 It follows that a bi-Lipschitz homeomorphism $f$ between two hyperbolic surfaces is also length-spectrum bounded, with length-spectrum constant bounded by any bi-Lipschitz constant of $f$.
        This implies that we have a natural inclusion
                \begin{equation}\label{eq:inclusion}
                \mathcal{T}_{bL}(H_0)\hookrightarrow \mathcal{T}_{ls}(H_0)
                \end{equation}
                and that  for any $x$ and $y$ of $\mathcal{T}_{bL}(H_0)$,  we have
                \begin{equation}\label{eq:lqb}
                d_{ls}(x,y)\leq d_{bL}(x,y).
                \end{equation}
                 In particular, the inclusion map in (\ref{eq:inclusion}) is  continuous.
        It will follow from the discussion in the next section that this inclusion is generally strict.

    Using the inclusion map in (\ref{eq:inclusion}),  the same example as in \ref{ex:ls}  shows that there exist pairs of marked hyperbolic structures on $X$ which are not bi-Lipschitz equivalent. Therefore the Teichm\"uller spaces $\mathcal{T}_{bL}(H_0)$ and $\mathcal{T}_{bL}(H_1)$ based at two different hyperbolic structures $H_0$ and $H_1$ on $X$ are in general setwise different.

    \begin{question} Characterize the hyperbolic surfaces $H_0$ on an infinite-type hyperbolic surface
 such that the inclusion
    in  (\ref{eq:inclusion}) is strict.
         \end{question}

\begin{definition} The {\it bi-Lipschitz mapping class group} of the surface $X$ relatively to the base hyperbolic metric $H_0$, denoted by $\mathrm{MCG}_{bL}(H_0)$, is the subgroup  of $\mathrm{MCG}(X)$  consisting of homotopy classes of self-homeomorphisms of $X$ that are bi-Lipschitz with respect to the metric $H_0$.
\end{definition}

From (\ref{eq:lqb}), we have an inclusion
\[ \mathrm{MCG}_{bL}(H_0)\subset \mathrm{MCG}_{ls}(H_0).
\]
Since we know that  $\mathrm{MCG}_{ls}(H_0)\varsubsetneq\mathrm{MCG}(X)$ (Proposition \ref{prop:ne}), we have 
 $\mathrm{MCG_{bl}}(H_0)\varsubsetneq\mathrm{MCG}(X)$.

The group $\mathrm{MCG}_{bL}(H_0)$ acts naturally on the Teichm\"uller space $\mathcal{T}_{bL}(H_0)$, by right-composition with the inverse at the level of markings (see Formula (\ref{eq:action-pre}) above), and this action preserves the metric   $d_{bL}$.

        Denoting as before  $\mathrm{MCG}_{f}(X)$ the subgroup of $\mathrm{MCG}(X)$ consisting of the homotopy classes of homeomorphisms which are supported on a subsurface of finite type of $X$, we have the following.

            \begin{proposition} There is a natural inclusion  $\mathrm{MCG}_{f}(X)\subset  \mathrm{MCG}_{bl}(H_0)$.
            \end{proposition}
            \begin{proof}
            As in the proof of Proposition \ref{prop:inclusion}, we use the fact that any element of the mapping class group of a surface of finite topological type is the product of a finite number of Dehn twist elements. A Dehn twist homeomorphism along a closed curve can be represented by a bi-Lipschitz homeomorphism supported on a compact annular neighborhood of that curve. The composition of these bi-Lipschitz homeomorphisms is a bi-Lipschitz homeomorphism of $X$.

            \end{proof}

        \section{The quasiconformal Teichm\"uller space of a surface of infinite type}\label{s11}
The Riemann surfaces, that is, the complex structures, that we shall consider will often be induced by hyperbolic metrics, and for this reason we shall use a metric definition of a quasiconformal map.
        We recall that if $f:S\to S'$ is a homeomorphism between two surfaces $S$ and $S'$ equipped with Riemannian metrics, and if $f$ is differentiable at a point $x$ of $S$, then the differential of $f$ at $x$ is an $\mathbb{R}$-linear map from the tangent space $T_xS$ of $S$ at $x$ to the tangent space $T_{f(x)}S'$ of $S'$ at $f(x)$, and therefore, such a linear map sends circles in $T_x S$ centered at the origin to ellipses in $T_{f(x)}S$ centered at the origin. The {\it quasiconformal dilatation $K(f)_x$ of $f$ at $x$} is then equal to the {\it eccentricity} of an image ellipse, that is, the ratio of the large axis to the small axis of this ellipse. This quantity is independent of the choice of the circle in $T_xS$ centered at the origin.
        Now if $f:S\to S'$ is a homeomorphism, then the {\it quasiconformal constant} of $f$ at a point $x$ in $S$ where $f$ is not necessarily differentiable is defined as follows.        For $r>0$, we consider the quantities
        \[L_f(x,r)=\sup\{d_{S'}(f(x),f(y))\ \vert \ d_S(x,y)=r\}\]
        and
        \[l_f(x,r)=\inf\{d_{S'}(f(x),f(y))\ \vert \ d_S(x,y)=r\}.\]
        The {\it quasiconformal dilatation} $K(f)_x$ of $f$ at $x$ is the element of $\mathbb{R}\cup\{\infty\}$ defined as
        \[ K(f)_x=\lim\sup_{r\to 0}\frac{L_f(x,r)}{l_f(x,r)}.\]

 In the case where $f$ is differentiable at $x$, the two definitions coincide.

        The quasiconformal dilatation of $f$ at $x$ depends only on the conformal classes of the Riemannian metrics on $S$ and $S'$. Therefore we can talk about the quasiconformal constant of a homeomorphism between Riemann surfaces, by choosing Riemannian metrics in the given conformal classes and using the above definition.

        The {\it maximal dilatation} $K(f)$ of a homeomorphism $f:S\to S'$ between two Riemann surfaces $S$ and $S'$, also called the {\it quasiconformal dilatation}, or, in short, the {\it dilatation} of $f$, is defined as the supremum of the quasiconformal dilatation $K(f)_x$ of $f$ over all points $x$ of $S$.
        The homeomorphism $f$ is said to be {\it quasiconformal} if its dilatation is finite.

       \begin{definition}
        Consider a Riemann surface structure $S_0$ on $X$. Its {\it  quasiconformal Teichm\"uller space}, $\mathcal{T}_{qc}(S_0)$, is defined as the set of equivalence classes $[f,S]$ of pairs $(f,S)$, where $S$ is a Riemann surface homeomorphic to $X$ and $f:S_0\to S$ a quasiconformal homeomorphism (called the marking of $S$), and where two pairs $(f,S)$ and $(f',S')$ are considered to be equivalent (and are said to be {\it conformally equivalent})  if there exists a conformal homeomorphism $f'':S\to S'$ homotopic to $f'\circ  f^{-1}$.
        \end{definition}
         The space $\mathcal{T}_{qc}(S_0)$ is equipped with the Teichm\"uller metric, defined as follows. Given two elements $[f,S]$ and $[f',S']$ of $\mathcal{T}_{qc}(S_0)$ represented by marked conformal surfaces $(f,S)$ and $(f',S')$,  their {\it quasiconformal distance} (also called the  {\it Teichm\"uller distance}, because it was introduced by Teichm\"uller) is defined as
        \begin{equation}\label{eq:qc} d_{qc}([f,S],[f',S'])=\frac{1}{2}\log \inf \{K(f'')\}
        \end{equation}
        where the infimum is taken over  maximal dilatations $K(f'')$ of homeomorphisms $f'':S\to S'$ homotopic to $f'\circ  f^{-1}$.

        The equivalence class of the marked Riemann surface $(\mathrm{Id},S_0)$ is the {\it basepoint} of $\mathcal{T}_{qc}(S_0)$. 
        
        It is important to note that the definition of the Teichm\"uller distance given in (\ref{eq:qc}) is the same as in the case of surfaces of finite conformal type, but that for surfaces of infinite type, choosing different basepoints may lead to different Teichm\"uller spaces. (But choosing a different basepoint by just applying a quasiconformal map leads to the same Teichm\"uller space, up to biholomorphism.) 
 Any homeomorphism between two Riemann surfaces of finite conformal type is homotopic to a differentiable map, which is automatically quasiconformal, but there exist homeomorphisms between surfaces of infinite conformal type that are not homotopic to quasiconformal maps. We shall discuss this fact in more detail below.

        We refer to Nag \cite{Nag} for an exposition of the quasiconformal theory of infinite-dimensional Teichm\"uller spaces.  In particular, it is known that the quasiconformal metric is complete.

      The quasiconformal dilatation and the bi-Lipschitz constant of a homeomorphism are related to each other. We shall see precise relations between these two quantities (see Theorem \ref{Th:Thurston} below).

    By the uniformization theorem, given any conformal structure $S_0$ on $X$, there is a hyperbolic structure $H_0$ on $X$ whose   underlying conformal structure is $S_0$.  We shall also denote the space $\mathcal{T}_{qc}(S_0)$ by  $\mathcal{T}_{qc}(H_0)$, meaning that we consider the conformal structure represented by the hyperbolic structure $H_0$.
 
From Wolpert's inequality (\ref{eq:Wolpert}) , there is a natural inclusion map
\begin{equation}\label{eq:inclusion2}
\mathcal{T}_{qc}(H_0)\hookrightarrow \mathcal{T}_{ls}(H_0).
\end{equation}

This inclusion map is continuous since Wolpert's result (\ref{eq:Wolpert}) also implies that for any two elements $x$ and $y$ in $\mathcal{T}_{qc}(S_0)$, we have
\begin{equation}\label{eq:lsqc}
d_{ls}(x,y)\leq d_{qc}(x,y).
\end{equation}

  \begin{proposition}
      In general, we have $\mathcal{T}_{ls}(H)\not=\mathcal{T}_{qc}(H)$.
        \end{proposition}

\begin{proof}
It suffices to construct a hyperbolic structure $H$ on $X$ and a self-homeomor\-phism of $(X,H)$ which is not quasiconformal and which is length-spectrum bounded. This is done in Proposition \ref{prop:suff} below.
\end{proof}

Concerning the comparison between the quasiconformal distance and the bi-Lipschitz distance, there is the following.
\begin{theorem} [Thurston \cite{Thurston1997} p. 268]\label{Th:Thurston} For every hyperbolic structure $H$, we have the set-theoretic equality
\begin{equation}\label{eq:T}
 \mathcal{T}_{qc}(H)=\mathcal{T}_{bL}(H)
\end{equation}
and  there exists a constant $C$ such that
for every $x$ and $y$ in $\mathcal{T}_{qc}(H)$, we have
\begin{equation}\label{eq:d}
d_{qc}(x,y)\leq d_{bL}(x,y)\leq C d_{qc}(x,y).
\end{equation}
\end{theorem}

It may be worth recalling that we are considering, as in Thurston's setting, complete hyperbolic structures.

The equality between the Teichm\"uller spaces in (\ref{eq:T}) follows from the inequalities between distances in (\ref{eq:d}).
The left hand side inequality in (\ref{eq:d}) follows from the fact that the quasiconformal dilatation is the supremum of an infinitesimal distortion at each point, and the bi-Lipschitz constant measures distortion of distances at all scales.  
The existence of a constant $C$ for which the right hand side inequality is satisfied is non-trivial, and it can be proved using the Douady-Earle quasiconformal barycentric extension theory  in \cite{DE} (see  \cite{Thurston1997} p. 268).

 \begin{corollary} For any hyperbolic structure $H$ on a surface $X$, we have the set-theoretic equality
$MCG_{qc}(H)=MCG_{bL}(H)$.
\end{corollary}

Douady and Earle gave in  \cite{DE} a proof of the fact that any quasiconformal Teichm\"uller space $\mathcal{T}_{qc}(H_0)$ is contractible (see \cite{DE} Theorem 3, where this result is also attributed to Tukia). We ask the following:

\begin{question}
Are the Teichm\"uller spaces $\mathcal{T}_{ls}(H_0)$ contractible ? Are they even connected ? 
\end{question}

Again, using the inclusion map in (\ref{eq:inclusion2}) and taking
the conformal structures $S_0$ and $S_1$ underlying the two
hyperbolic structures $H_0$ and $H_1$ of
    Example \ref{ex:ls}, we obtain an element $H_1\not\in \mathcal{T}_{qc}(H_0)$. This gives an example of two Teichm\"uller spaces $\mathcal{T}_{qc}(S_0)$ and $\mathcal{T}_{qc}(S_1)$ whose basepoints are two homeomorphic Riemann surfaces $S_0$ and $S_1$, such that 
     $\mathcal{T}_{qc}(S_0)\not=\mathcal{T}_{qc}(S_1)$ setwise.

\begin{definition} Let $S$ be a Riemann surface structure on the topological surface $X$. The {\it quasiconformal mapping class group} of $S$, denoted by $\mathrm{MCG}_{qc}(S)$, is the subgroup  of $\mathrm{MCG}(X)$  consisting of homotopy classes of homeomorphisms that are quasiconformal with respect the structure $S$.
 \end{definition}

The group $\mathrm{MCG}_{qc}(S)$ naturally acts on the Teichm\"uller space $\mathcal{T}_{qc}(S)$, by right-composition with the inverse at the level of markings (see the action described in (\ref{eq:action-pre})), and this action preserves the metric $d_{qc}$.

The homomorphism of the group $\mathrm{MCG}_{qc}(S)$ into the group of isometries of $\mathcal{T}_{qc}(S)$ is injective (see \cite{Epstein2000} and \cite{EGL}).

A theorem whose most general version is due to V. Markovic  (see \cite{M}) asserts that any isometry of the space $\mathcal{T}_{qc}(S)$ is induced by an element of   $\mathrm{MCG}_{qc}(S)$. Special cases of this result had been obtained by several people, including
Earle and Kra \cite{EK},  Lakic \cite{L} and Matsuzaki \cite{Mat}. The earliest version of this theorem, in the case of finite-type surfaces, is due to Royden \cite{Ro}.

\begin{proposition}\label{prop:suff} Let $H$ be a hyperbolic structure on $X$ such that
there exists a sequence of homotopy classes of disjoint essential
closed curves on $X$ whose length tends to 0.
Then,
 $\mathrm{MCG}_{qc}(H)\varsubsetneq\mathrm{MCG}_{ls}(H)$.
  \end{proposition}

\begin{proof}
We produce an element $T\in \mathrm{MCG}_{ls}(H)$ with $T\not\in\mathrm{MCG}_{qc}(H)$.

Without loss of generality, we can assume that there exists a sequence
 $\alpha_n, n=1,2 \ldots$, of homotopy classes of disjoint essential
closed curves on $X$
 whose lengths satisfy  $l_H(\alpha_n)=\varepsilon_n$ with $e^{-(n+1)^2}<\varepsilon_n<e^{-n^2}$.

Let $$ t_n= \left[\frac{\log |\log
\epsilon_n |}{\epsilon_n}\right]+1, \ n=1, 2, ... $$ ($[x]$ denoting the
integral part of the real number $x$).

 For each $n=1, 2, ...$, let $\tau_n$ be the $t_n$-th power of the positive Dehn
twist about $\alpha_n$. We take all the positive Dehn twists to be supported on disjoint annuli, and
we define $T\in\mathrm{MCG}(X)$ to be the infinite composition
$T=\tau_1 \circ \tau_2 \circ \ldots$.  From \cite{Liu2008},
we have that $T$ is not in
$MCG_{qc}(X)$.

Now we show that $T \in
MCG_{ls}(X)$. We need to show that
\[\label{eq:lss}\log \sup_{\alpha\in\mathcal{S}(X)} \left\{
\frac{l_{H}(f(\alpha))}{l_{H}(\alpha)},
\frac{l_{H}(\alpha)}{l_{H}(f(\alpha))}\right\}<\infty.\]
The proof parallels the proof of Proposition \ref{prop:t0}.

 Let $\alpha$
be an arbitrary homotopy class of essential curves in $X.$

If $i(\alpha, \alpha_n)=0, n=1,2,
\ldots$, then $\alpha=T(\alpha)$ and  $l_{H} (T (\alpha))= l_{H}
(\alpha)$.

Suppose that $i(\alpha, \alpha_n) \neq 0$
for some $n$. Without loss of generality, we can assume that
$i(\alpha, \alpha_n) \neq 0$, for all $n=1,2, \ldots,N$. By the  Collar Lemma, on any hyperbolic surface, any closed geodesic whose length $\epsilon$ is sufficiently small has an embedded collar neighborhood of width $\vert \log\epsilon\vert$. Thus, we can write
$$l_{H}(\alpha)= \Sigma_{n=1}^N i(\alpha, \alpha_n) |\log\varepsilon_n|+ B_{\alpha}$$
where $ B_{\alpha}>0$ is a constant that depends on the closed curve $\alpha$.

From the definition of a Dehn twist, we have
\[\displaystyle \Sigma_{n=1}^N i(\alpha, \alpha_n) |\log\varepsilon_n| + B_{\alpha} -
\Sigma_{n=1}^N i(\alpha, \alpha_n) t_n \varepsilon_n \leq
l_{H}(T(\alpha))\]
\[\displaystyle \leq \Sigma_{n=1}^N i(\alpha, \alpha_n)
|\log\varepsilon_n| + B_{\alpha} + \Sigma_{n=1}^N i(\alpha,
\alpha_n) t_n \varepsilon_n .\]

  Therefore we have
  \begin{eqnarray*}
\frac{l_{H}(T(\alpha))}{l_{H}(\alpha)}&\leq&
 \frac{\Sigma_{n=1}^N i(\alpha, \alpha_n) |\log\varepsilon_n| +
B_{\alpha} + \Sigma_{n=1}^N i(\alpha, \alpha_n) t_n \varepsilon_n}
{\Sigma_{n=1}^N i(\alpha, \alpha_n) |\log\varepsilon_n|+
B_{\alpha}}
 \\&=&
  1+ \frac{\Sigma_{n=1}^N i(\alpha, \alpha_n) t_n \varepsilon_n}{\Sigma_{n=1}^N i(\alpha, \alpha_n) |\log\varepsilon_n|+ B_{\alpha}}
  \\&\leq&
   1+ \frac{\Sigma_{n=1}^N i(\alpha, \alpha_n)\log|\log
\varepsilon_n|}{\Sigma_{n=1}^N i(\alpha, \alpha_n) |\log
\varepsilon_n|}
  \\&\leq&
 1+ {\Sigma_{n=1}^N} \frac {\log|\log \varepsilon_n|} {|\log
\varepsilon_n|}
  \\&\leq& 1+ 2 \Sigma_{n=1}^N \frac{\log (n+1)}{n^2}
  \\&\leq & 1+ 2 \Sigma_{n=1}^{\infty} \frac{\log (n+1)}{n^2}<\infty
\end{eqnarray*}
which is bounded independently of $\alpha$ and $n$. In the same way, we can prove that
$\frac{l_{H}(\alpha)}{l_{H}(T(\alpha))}$ is also bounded
independently of $\alpha$ and $n$. This implies that $T \in
MCG_{ls}(H)$.

\end{proof}

        \begin{question}    Proposition \ref{prop:suff} gives a sufficient condition  under which the inclusion in (\ref{eq:inclusion2}) is strict.  The question is to find  sufficient and necessary conditions for the same fact.
        \end{question}

Recall that $\mathrm{MCG}_f(X)$ denotes the group of homotopy classes of homeomorphisms of $X$ which are supported on subsurfaces of finite analytical type.

            \begin{proposition}
            For any conformal structure $S$ on $X$, we have a natural inclusion  $\mathrm{MCG}_{f}(X)\subset  \mathrm{MCG}_{qc}(S)$.
            \end{proposition}

            \begin{proof} Any homeomorphism of a subsurface of $X$ of finite analytical type is homotopic to a diffeomorphism, and therefore to a quasiconformal homeomorphism. This quasiconformal homeomorphism can be extended to a quasiconformal homeomorphism of $X$ with the same quasiconformal dilatation constant.
            \end{proof}

 It was shown in \cite{Liu2008} that the
following is a necessary condition for the equality
$\mathcal{T}_{qc}(H_0)=\mathcal{T}_{ls}(H_0)=\mathcal{T}_{bL}(H_0)$: there exists a positive constant $\epsilon$, such that
for any element $\alpha$ in $\mathcal{S}(X)$, $l_{H_0}(\alpha)\ge \epsilon$.

Let $C_i$ ($i=1,2,\ldots$) be the collection of the simple closed curves that are boundaries of pairs of pants in the pants decomposition $\mathcal{P}$ of the surface $X$. We recall that Condition (\ref{LM}) above,  for a hyperbolic structure $H$ on $X$, states the following:
\[
 \displaystyle \exists M>0 , \forall i=1,2,\ldots,
 \frac{1}{M}\leq l_{H}(C_i)\leq M.
\]
   
    \begin{theorem}\label{th:Sh}
    Suppose that the base hyperbolic metric $H$ satisfies Condition (\ref{LM}) and
    let $x_n, n=0,1, \ldots, $ be a sequence of points in $\mathcal{T}_{qc}(H)$. Then the following are
    equivalent:
\begin{enumerate}[(a)]
  \item\label{a} $\lim_{n \to \infty}d_{qc}(x_n, x_0)=0$;
 \item \label{b} $\lim_{n \to \infty}d_{ls}(x_n, x_0)=0$;
    \item \label{c}  $\lim_{n \to \infty}d_{bL}(x_n, x_0)=0$.
\end{enumerate}
    \end{theorem}

    \begin{proof} From (\ref{eq:d}), we have (\ref{a})$\Rightarrow$(\ref{b}).
    From Theorem \ref{Th:Thurston}, we have (\ref{a})$\iff$(\ref{c}). Finally, from Shiga's result in \cite{Shiga2003} (proof of Theorem 1.2), we have (\ref{b})$\Rightarrow$(\ref{a}).
     \end{proof}

 We consider for a while surfaces of finite type.

The following follows from Lemma \ref{lem:Liu}.

\begin{lemma} \label{lem:ccc} Let $S$ be a  finite type Riemann surface
of genus $g$ with $n$ punctures and $c$ boundary components, and let $\epsilon >0$. Suppose that  there exists $n_0>0$ such that $g+n+c \le n_0$. For any $x$ and $y$ in $\mathcal{T}^{\epsilon}(S)$, there  exists  a constant
$N(n_0, \epsilon)$, depending only on $n_0$ and $\epsilon$, such that
\[
 d_{ls}(x,y) \le d_{qc}(x,y)
 \le 4 d_{ls}(x,y)+N(n_0, \epsilon).
\]
 \end{lemma}

 \begin{proof} For  $c=0$, this lemma is a reformulation of Lemma \ref{lem:Liu}. For $c\not=0$, the result follows by taking doubles of surfaces along geodesic boundary components.
 \end{proof}

Let $S$ be a Riemann surface equipped with a hyperbolic metric and let $\alpha$ be a simple closed geodesic on $S$. We recall that a time-$t$ Fenchel-Nielsen deformation of $S$ with respect to $\alpha$ is a Riemann surface $S_t$ obtained by cutting $S$ along $\alpha$ and regluing the cut off borders after twisting by  a hyperbolic distance $t$ along $\alpha$. If $S$ is a marked surface, then $S_t$ is also equipped with a natural marking. We refer to \cite{IT} for background material.

The following two lemmas are inspired from Lemma 3.1 in \cite{Shiga2003}, and a corresponding lemma in \cite{Liu2008}.

\begin{lemma} \label{lem:w1} Let $S$ be a hyperbolic surface and let $\alpha$ be a simple
closed curve on $S$ which is not homotopic to a boundary component
of $S$. Let $S_t$ be the time-$t$
Fenchel-Nielsen deformation of $S$ about $\alpha$, let $\tau^t_{\alpha}: S \to S_t$ be a quasiconformal homeomorphism and let $K(t)$ be its
dilatation. Then $\lim_{t \to \infty} K(t)=
\infty$.
\end{lemma}

\begin{proof} Let $\rho$ be the hyperbolic distance on $S$. For $s>0$, let
\[
W_{\alpha}=\{p \in S \vert \rho (p, \alpha)<s\}.
\]
The positive constant $s$ is chosen to be so small that $W_{\alpha}$
is a tubular neighborhood of $\alpha$ in $S$. We construct a quasiconformal homeomorphism $g_t$ from $S$ to $S_t$ which represents the positive twisting along $\alpha$ by $t$. This map is best described by looking to a lift to the universal covers of the two surfaces $S$ and $S_t$.

  We identify the universal covers of $S$ and $S_t$ with the upper-half model of the hyperbolic plane, $\mathbb{H}^2$. We denote by $\pi: \mathbb{H}^2
\to S$ the covering map. We consider connected components of ${\pi}^{-1} (W_{\alpha})$ and of ${\pi}^{-1} (\alpha)$ in $\mathbb{H}^2$.   Without loss of generality, we may assume that the
connected component of ${\pi}^{-1} (\alpha)$ is the positive half of
the imaginary axis. Let $\tilde W_{\alpha}$ be the connected
component of ${\pi}^{-1} (W_{\alpha})$ containing the positive half
of the imaginary axis. For a suitable $\theta$ with $0 < \theta
< \frac{\pi} {2}$, we can express $\tilde{W_{\alpha}}$ as
$$
\tilde{W_{\alpha}}=\{z \in \mathbb{H}^2 \ \vert \  \frac {\pi} {2} -\theta< \arg z <
\frac{\pi} {2} +\theta \}.
$$

Let $\omega_t$ be the lift of a  quasiconformal mapping in
the homotopy class of $g_t,$ normalized by $\omega_t (0)=0,
\omega_t(i)=i$ and $\omega_t (\infty)= \infty.$ It is well known that
$\omega_t$ can be extended to a homeomorphism of $\overline{\mathbb{H}^2}$ and
that the boundary mapping $\omega_t \vert_{\mathbb{R}}$ depends only on the
homotopy class of $g_t$.

 Finally, we let $\delta$ be the geodesic in $\mathbb{H}^2$ which (as a set) is defined as
 \[\delta =\{ z\in \mathbb{H}^2, \vert z\vert = 1.\}\]

 Let $z_1$ and $z_2$,  $(\Re {z_1} < 0 < \Re {z_2})$ be the points of
$\delta \cap {\tilde W_{\alpha}}$. From the fact that $g_t$ is the positive time-$t$
Fenchel-Nielsen deformation, we have $\omega_t (z_1)=z_1$, 
$\omega_t (z_2)=e^{t} z_2,$ and $\omega_t (\delta \cap {\tilde
W_{\alpha}})$ is an arc connecting $z_1$ and $e^{t}z_2$ in the
 component $\tilde W_{\alpha}$ (see Figure \ref{arc}). 
 There is a similar description of the image by $\omega_t$ of a
 subarc of $\delta$ in each component of ${\pi}^{-1}
(W_{\alpha})$. Let $z_1=-\sin \theta + i \cos \theta$ and $z_2=\sin \theta
+ i \cos \theta$. Then,  
\[\omega_t (z_2)=e^t (\sin \theta + i \cos
\theta)\]
and  
\[\Re \omega_t(z_2)-\Re z_2=(e^t-1) \sin\theta.\]

  \begin{figure}[!hbp]
\centering
\psfrag{2}{\small $-1$}
\psfrag{0}{\small $0$}
\psfrag{1}{\small $1$}
\psfrag{O}{\small $O$}
\psfrag{w}{\small $z_1$}
\psfrag{k}{\small $e^tz_2$}
\psfrag{e}{\small $z_2$}
\psfrag{t}{\small $\theta$}
\includegraphics[width=.60\linewidth]{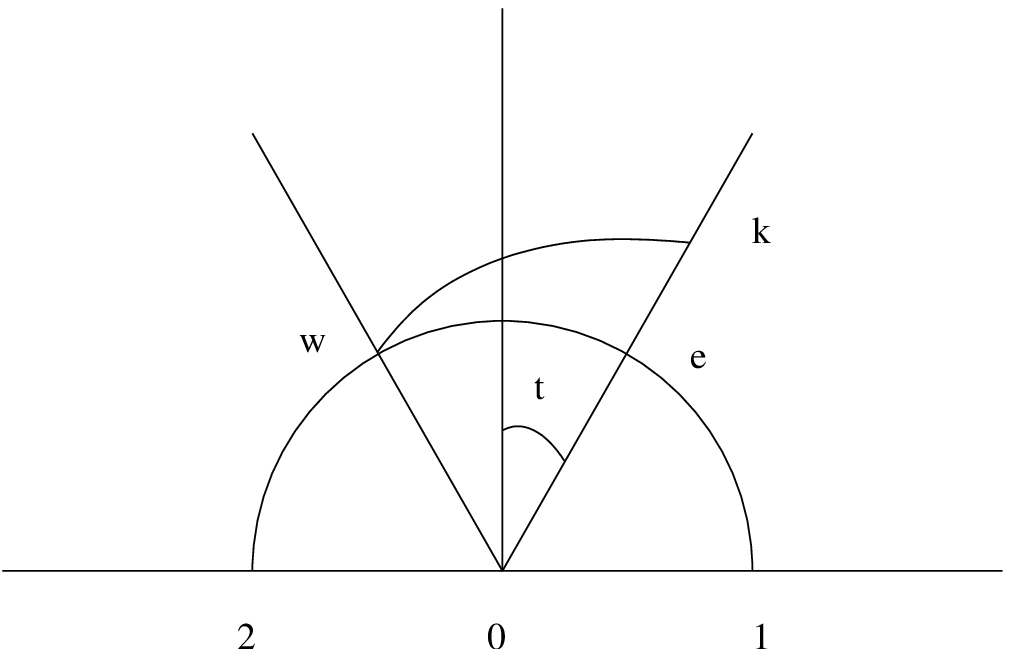}
\caption{\small {}}
\label{arc}
\end{figure}

Now we study the action of the quasiconformal map $\omega_t$ on the line at infinity $\mathbb{R}$. From
the above construction we have 
\[\omega_t (1)-1>\Re \omega_t(z_2)-\Re
z_2=(e^t-1) \sin\theta\]
and 
\[-1<\omega_t(-1)<0.\]
Summing up, we have

\begin{equation}\label{eq:w}
-1 < \omega_t (-1) < 0 <1+(e^t-1) \sin\theta < \omega_t (1).
\end{equation}
 We use the following definition of the cross ratio of quadruples in $\mathbb{H}^2$:
$$[a, b, c, d]= \frac{(a-b)(c-d)}{(a-d)(c-b)}.
$$
This definition extends to a notion of cross ratio defined for quadruples of points at infinity, and we have
$$
[-1, 0, 1, \infty]= -1.
$$

We shall use the geometric definition of a quasiconformal mapping, that involves the ratios of moduli of quadrilaterals, and we refer to Ahlfors' book \cite{Ahlfors} for this definition.

We can talk about the conformal modulus of ideal quadrilaterals (quadrilaterals having their vertices at infinity), and we have the following formula:
\begin{equation}\label{eq:w1}
\mathrm{Mod} (H (-1, 0, 1, \infty))=1.
\end{equation}
We set
$$\nu_t = [ \omega_t(-1), \omega_t (0), \omega_t (1), \omega_t (\infty)]=[\omega_t (-1), 0, \omega_t (1),
\infty]= \frac{\omega_t (-1)}{\omega_t (1)}.
$$
Then
$$|\nu_t|= |\frac{\omega_t (-1)}{\omega_t (1)}| \leq \frac{1}{1+(e^t-1) \sin\theta} \to 0, \ t \to
\infty.
$$
Therefore, we have the following.
 \begin{equation}\label{eq:w2}
 \mathrm{Mod} (H (\omega_t
(-1), \omega_t (0), \omega_t (1), \omega_t (\infty))) \to 0, \ as \
t \to \infty.
\end{equation}
By the geometric definition of quasiconformal mappings, we get
\begin{equation}\label{eq:w3}
\frac{1}{K(\omega_t)} \leq \frac{\mathrm{Mod} (H (\omega_t (-1), \omega_t (0), \omega_t
(1), \omega_t (\infty)))}{\mathrm{Mod} (H (-1, 0, 1, \infty))}.
\end{equation} 
Summing up, (\ref{eq:w1}), (\ref{eq:w2}) and (\ref{eq:w3}) give
\begin{equation}\label{eq:w4}
K(\omega_t) \to \infty, \ n
\to \infty. \end{equation}

Since any quasiconformal mapping in the homotopy class of $\omega_t$ induces the same boundary map as $\omega_t$, its quasiconformal dilatation $K(t)$ tends to infinity as $t\to \infty$.

\end{proof}

\begin{lemma} \label{lem:w2} Let $S$ be a Riemann surface and let $\alpha_n$ be a
sequence of simple closed curves on $S$ which are not homotopic to a
boundary component of $S$. Given an unbounded sequence $t_n$ of real numbers, let $\tau^{t_n}_{\alpha_n}: S \to
S_{t_n}$ be the time-$t_n$ Fenchel-Nielsen deformation about
$\alpha_n$, and let $K(t_n)$ be the dilatation of
$\tau^{t_n}_{\alpha_n}$. Then $\lim_{t_n \to \infty} K(t_n)=
\infty$.
\end{lemma}

 Lemma \ref{lem:w2} can be proved  in the same way as Lemma \ref{lem:w1}.

\begin{theorem}
    Suppose that the base hyperbolic metric $H_0$ satisfies Condition (\ref{LM}).
    Then, $\mathcal{T}_{qc}(H_0)=\mathcal{T}_{ls}(H_0)=\mathcal{T}_{bL}(H_0).$

    \end{theorem}

\begin{proof}

    It is easy to see that if $H_0$ satisfies Condition (\ref{LM}), then every element of $\mathcal{T}_{ls}(H_0)$ satisfies that condition (but of course, not necessarily with the same constant $M$). Since the hyperbolic length of the image of a homotopy class $\alpha$ of closed curves by a $K$-Lipschitz homeomorphism is bounded from above by $K$ times the hyperbolic length of $\alpha$, it follows that  if $H_0$ satisfies Condition (\ref{LM}), then every element of $\mathcal{T}_{bL}(H_0)$ satisfies that condition.  Finally, using Wolpert's inequality   (\ref{eq:Wolpert}), it is also easy to show that if $H_0$ satisfies Condition (\ref{LM}), then every element of $\mathcal{T}_{qc}(H_0)$ satisfies the same condition. (One can also use the set-theoretic equality $\mathcal{T}_{bL}(H_0)=\mathcal{T}_{qc}(H_0)$ and their inclusion in $\mathcal{T}_{ls}(H_0)$.)

    From  (\ref{eq:inclusion}), we have $\mathcal{T}_{bL}(H_0)\subset \mathcal{T}_{ls}(H_0)$. From (\ref{eq:inclusion2}), we have $\mathcal{T}_{qc}(H_0)\subset \mathcal{T}_{ls}(H_0)$. From Theorem  \ref{Th:Thurston}, we have $\mathcal{T}_{bL}(H_0)= \mathcal{T}_{qc}(H_0)$.
    Therefore, it remains to show that
 $\mathcal{T}_{ls}(H_0)\subset \mathcal{T}_{qc}(H_0)$. We now prove this by showing that any two hyperbolic structures at bounded length spectrum distance are at bounded quasiconformal distance.

We consider two decompositions of $X$ by surfaces with disjoint interiors, 
 $X=\bigcup_{n=1}^{\infty} X_n$ and  $X=\bigcup_{n=1}^{\infty} X'_n$, satisfying the following properties:
 
\begin{itemize}  \item For each $n=1,2, \ldots$, each of the surfaces $X_n$ and
$X'_n$ is of finite type, and  the sum of its  genus, number of punctures and 
number of boundary components is bounded by an integer $n_0$ that is independent of $n$.
\item If $S_i, i=1, 2, \ldots$  and $S'_i, i=1, 2, \ldots$ are the
decomposition curves of the surface decompositions  $X=\bigcup_{n=1}^{\infty} X_n$ and  $X=\bigcup_{n=1}^{\infty} X'_n$ respectively, then
 $\{S_i\}_{i=1}^{\infty} \subset \{C_i\}_{i=1}^{\infty}$ 
 and  $\{S'_i\}_{i=1}^{\infty} \subset \{C_i\}_{i=1}^{\infty}$,
where $C_i$, $i=1,2,\ldots$ are the curves in the pants
decomposition $\mathcal{P}$ of the surface $X$ satisfying condition
 (\ref{LM}).
 
 \item 
We have $\{ \bigcup_{i=1}^{\infty} S_i
\} \bigcap \{ \bigcup_{i=1}^{\infty} S'_i \}= \emptyset$.
\end{itemize}

 Let $x_1=[f_1, H_1]$ and $x_2=[f_2, H_2]$ be two elements of the Teichm\"uller space 
 $\mathcal{T}_{ls}(H_0)$ and let $M_0=d_{ls}(x_1, x_2)$

 For $i=0,1,2$ and for each $n=1,2,\ldots$, let $H_i^n$, be the restriction of the hyperbolic
structure $H_i$, on $X_n$, let $f_i^n :
H_0^n \to H_i^n$, be the restriction of $f_i$ to $H_0^n$ and let
$x_i^n=[f_i^n, H_i^n]$.

For any $i=1,2$ and for each $n=1,2,\ldots$, we have $x_i^n \in \mathcal{T}^{\epsilon}(H_0^n)$, where
$\epsilon$ only depends on the the constant $M$ that appears in Condition  (\ref{LM}). From Lemma \ref{lem:ccc}, we can write, for any $n=1,2,\ldots$, 
\[
 d_{ls}(x_1^n, x_2^n) \le d_{qc}(x_1^n, x_2^n)
 \le 4 d_{ls}(x_1^n, x_2^n)+N(n_0, \epsilon)
\le 4 M_0+N(n_0, \epsilon).
\]

It remains to show that the amounts of twist along the curves $S_i, i=1, \ldots$, are
bounded. 

For an arbitrary integer $p\geq 1$, we can find an integer $n_p$ such that $S_{p} \subset
X'_{n_{p}}$.  For $i=0,1,2$, let $H_{i_{p}}'$, be the restriction of $H_i$ to  $X'_{n_{p}}$ and  let $f'_{i_{p}}: H'_{0_{p}} \to H'_{i_{p}}$, $i=0,1,2$, be the map obtained by restriction of $f$.
Let $x'_{i_{p}}=[f'_{i_{p}}, H'_{i_{p}}]$, $i=0,1,2$, Then, by the same arguments that we used above, we have
\begin{equation}
\label{ineq}
d_{qc}(x'_{1_{p}}, x'_{2_{p}}) \le 4M_0+N(n_0, \epsilon).
\end{equation}
From Lemmas   \ref{lem:w2} and \ref{lem:w1},  the amount of twist about
$C_{i_{0}}$, $i=0,1,2$ is bounded by a constant $B$, depending only on $M, M',
n_0$, since  otherwise, if this amount of twist had no bounds, Lemma
\ref{lem:w1}  would imply that the dilatation of the corresponding quasiconformal
mapping tends to infinity. This contradicts Inequality (\ref{ineq}).

Since for an arbitrary integer $p\geq 1$  
$d_{qc}(x'_{1p}, x'_{2p}) \le 4M_0 +N(n_0, \epsilon)$, and since the twists about the  boundary curves have a common upper bound,
 we have $d_{qc}(x_1, x_2) < \infty$.
     \end{proof}

\begin{theorem}
    Suppose that the base hyperbolic metric $H_0$ satisfies Condition (\ref{LM}).
    Let $x_n, n=0,1, \ldots, $ be a sequence of marked Riemann surface relative to
    the base hyperbolic structure $H_0$.  Then the following are
    equivalent:

   \begin{enumerate}[(a)]  
   \item \label{e1} $\lim_{n \to \infty}d_{qc}(x_n, x_0)=\infty$,

   \item \label{e2} $\lim_{n \to \infty}d_{ls}(x_n, x_0)=\infty$,

 \item  \label{e3}  $\lim_{n \to \infty}d_{bL}(x_n, x_0)=\infty$.
\end{enumerate}
    \end{theorem}

    \begin{proof} From (\ref{eq:d}) we know that (\ref{e1})  is equivalent to (\ref{e3}).
    From (9) we know that (\ref{e1})  implies (\ref{e2}). By  quoting Shiga's result as we did in the proof of Theorem \ref{th:Sh}, we have (\ref{e2})  implies (\ref{e1}).
     \end{proof}

Combining Theorems \ref{th:Sh} and \ref{th:co}, we obtain
     \begin{theorem}
    Suppose that the base hyperbolic metric $H$ satisfies Condition (\ref{LM}).
    Then $\mathcal{T}_{qc}(H)$, $\mathcal{T}_{ls}(H)$ and
    $\mathcal{T}_{bL}(H)$ are all complete.

    \end{theorem}

\begin{corollary} If the base  hyperbolic metric $H$
satisfies Condition (\ref{LM}), then $MCG_{qc}(H)=MCG_{ls}(H)$.
\end{corollary}

\end{document}